\documentclass[preprint,3p,12pt]{elsarticle}
\usepackage[utf8]{inputenc}
\usepackage{textcomp}

\usepackage{graphicx}
\usepackage{amsmath,amssymb}
\usepackage[version=4]{mhchem}
\usepackage{siunitx}
\usepackage{longtable,tabularx}
\usepackage[symbol]{footmisc}

\usepackage{bm}
\usepackage{amsfonts}
\usepackage{amscd}
\usepackage{hyperref}
\usepackage{pstricks}
\usepackage{subcaption}
\usepackage{framed,fancybox}
\usepackage{bbm}
\usepackage{multirow}
\usepackage{threeparttable}
\usepackage{rotating}
\usepackage{stfloats}
\usepackage{color,soul}
\usepackage{float}
\usepackage{afterpage}
\usepackage{listings}
\usepackage{xcolor}

\usepackage{amsthm}

\definecolor{codegreen}{rgb}{0,0.6,0}
\definecolor{codegray}{rgb}{0.5,0.5,0.5}

\lstdefinestyle{mystyle}{
  commentstyle=\color{codegreen},
  numberstyle=\tiny\color{codegray},
  keywordstyle=\color{blue},
  basicstyle=\ttfamily\small,
  breakatwhitespace=false,         
  breaklines=true,                 
  keepspaces=true,                 
  numbers=left,                    
  numbersep=5pt,                  
  showspaces=false,                
  showstringspaces=false,
  showtabs=false,                  
  tabsize=2
}

\allowdisplaybreaks

\journal{Journal of the Franklin Institute}

\begin{document}

\begin{frontmatter}

\title{{\bf An Adaptive Method for Optimal Control Problems} \\{\bf Constrained by Parabolic Differential Equations}}

\author[ufaffil]{Alexander M.~Davies\fnref{label1}}
\ead{alexanderdavies@ufl.edu}
\fntext[label1]{Ph.D. Candidate, Department of Mechanical and Aerospace Engineering.}
\affiliation[ufaffil]{organization={University of Florida},
            city={Gainesville},
            postcode={32611},
            state={FL},
            country={USA}}

\affiliation[afrlaffil]{organization={Air Force Research Laboratory},
            city={Eglin AFB},
            postcode={32542},
            state={FL},
            country={USA}}
            
\author[ufaffil]{Sara Pollock\fnref{label2}}
\ead{s.pollock@ufl.edu}
\fntext[label2]{Associated Professor, Department of Mathematics.}

\author[afrlaffil]{Miriam E.~Dennis\fnref{label3}}
\ead{miriam.dennis.1@us.af.mil}
\fntext[label3]{Research Engineer, Munitions Directorate.}

\author[ufaffil]{Anil V. Rao\corref{cor1}\fnref{label4}}
\ead{anilvrao@ufl.edu}
\cortext[cor1]{Corresponding author.}
\fntext[label4]{Professor, Department of Mechanical and Aerospace Engineering. AIAA Associate Fellow. AAS Fellow.}

\begin{abstract}
An adaptive direct collocation method is developed for solving optimal control problems constrained by parabolic partial differential equations. The partial differential equation is first reformulated in a variational setting, where the spatial domain is discretized using the $hp-$Galerkin finite element method. To address nonlinearities in the variational form, a ``Kirchhoff-like'' integral transformation is applied to linearize the dynamics. In the temporal dimension, an orthogonal collocation scheme, the $hp-$flipped Legendre–Gauss–Radau method, is employed to fully discretize the problem, yielding a large, sparse nonlinear programming problem. Upon solving the nonlinear programming problem, solution accuracy is assessed through an implicit residual estimation procedure. This approach evaluates the local error by solving auxiliary residual problems over selected subdomains, providing a novel means of error estimation within an orthogonal collocation framework for optimal control. Based on the computed error estimate, the mesh is adaptively refined or coarsened to meet a prescribed error tolerance. Mesh refinement is guided by the estimated regularity of the solution which is determined via the decay rate of the coefficients of a Legendre polynomial expansion. In overcollocated regions, a mesh reduction strategy is adapted from orthogonal collocation methods for application within the finite element framework. Numerical examples demonstrate that the proposed method can reduce the $L_2-$error by up to five orders of magnitude in both spatial and temporal dimensions.
\end{abstract}


\end{frontmatter}
\newpage

\section{Introduction}
Numerical methods for optimal control (OC) have primarily been developed on optimal control problems (OCPs) constrained by ordinary differential equations (ODEs). In recent decades, as the field and computational resources continue to grow, the study of numerical methods for OCPs constrained by partial differential equations (PDEs) has garnered more interest \cite{Lions1971,Troltzsch2024,ManzoniQuarteroni2021, BeckerKapp2000}. Some fields in which PDE-constrained OC have been applied include, but are not limited to: fluid mechanics \cite{Casas1998,Ghattas1997,AllahverdiPozo2016,Sritharan1998}, thermodynamics \cite{DebusscheFuhrman2007,CasasTroltzsch2020,CleverLang2012}, finance \cite{KafashDelavarkhalafi2016}, and epidemiology \cite{WangZhang2023}. In general, PDE-constrained OCPs are complex and generally lack analytical solutions due to the increased dimensionality and dynamical complexity of the underlying system. As a result, in most cases, efficient and accurate numerical methods are required to produce solutions to PDE-constrained OCPs.


In general, numerical methods for PDE-constrained OC fall under two categories: indirect and direct methods. In an indirect method, the OCP is transcribed into a Hamiltonian boundary value problem (HBVP) through the derivation of the necessary first-order optimality conditions from variational calculus. The HBVP can then be solved by techniques such as shooting \cite{HesseKanschat2009, CarraroGeiger2014,CarraroGeiger2015} or collocation \cite{MeidnerVexler2008a, BeckerBraack2007, SabehShamsi2016, MohammadizadehTehrani2019}. In a direct method, the OCP is transformed into a large, sparse nonlinear programming problem (NLP) through control and/or state parameterization. The NLP can then be solved by well-developed software such as \textit{SNOPT} \cite{GillMurray2002}, \textit{IPOPT} \cite{BieglerZavala2008}, or \textit{KNITRO} \cite{ByrdNocedal2007}.

Direct methods for PDE-constrained OC have grown increasingly popular in recent decades. The use of a direct method for parabolic PDE-constrained OC requires the use of suitable spatial and temporal discretizations. Some popular discretizations for PDE-constrained OC that have been used to discretize the spatial dimension include the finite difference method \cite{BuskensGriesse2006,Betts2020,KameswaranBiegler2008,FangVandewalle2022}, the finite element method \cite{KupferSachs1992, Heinkenschloss1996, ChenZhou2018, DaviesPollock2025}, the finite volume method \cite{HolmqvistMagnusson2016}, and orthogonal collocation methods \cite{JieZhu2023,KhaksarShamsi2017,AliShamsi2019}. In the temporal dimension, typical choices include Euler and Runge-Kutta methods \cite{FangVandewalle2022,ChenZhou2018, Betts2020}, spectral methods \cite{Nemati2018}, and orthogonal collocation methods \cite{KhaksarShamsi2017,JieZhu2023, DaviesPollock2026}. Following the choice of a suitable discretization, one may solve the resulting NLP via a control parameterization approach such as shooting and multiple shooting \cite{BuskensGriesse2006, FangVandewalle2022, ChenZhou2018} or by a state and control parameterization approach known as collocation \cite{Betts2020,KameswaranBiegler2008, Heinkenschloss1996,KupferSachs1992,Nemati2018, DaviesPollock2025,KhaksarShamsi2017, HolmqvistMagnusson2016, JieZhu2023, AliShamsi2019,DaviesPollock2026}.

Direct collocation is arguably the most powerful method for solving general OCPs \cite{Rao2009}. In a direct collocation method, the domain is partitioned into a mesh comprised of mesh intervals (temporal) and elements (spatial). Within each interval or element, the state is approximated by a set of basis or trial functions supported at a discrete set of points, and constraints are enforced at a specific subset of the support points called \textit{collocation points}. 

Historically, direct collocation methods have been formulated as $h-$methods or $p-$methods. In an $h-$method, the order of the state approximation within each interval/element is fixed, and convergence is achieved by increasing the number of intervals/elements and/or improving the placement of mesh points within the domain. In a $p-$method, the number of intervals/elements is fixed (typically one), and convergence is achieved by increasing the order of the state approximation within each interval/element \cite{Betts2020}. While both $h-$methods and $p-$methods have been used extensively and successfully, both approaches can have limitations in achieving high-accuracy solutions. With an $h-$method, the number of required elements may grow intractably large, and similarly, in a $p-$method, the order of the approximation within an element may grow unreasonably large. To mitigate the limitations of both $h-$ and $p-$methods, so-called $hp-$methods have been developed. In an $hp-$method, both the number and placement of the mesh points \textit{and} the order of the approximation within each element can be varied to achieve convergence. Classes of $hp-$methods have grown increasingly popular due to their robustness and compuational performance in comparison to $h-$ and $p-$methods alone. 

Two methods that can be cast as direct $hp-$collocation methods are the finite element method and the orthogonal collocation method. In either case, both the degree of the state approximation and the number of elements/intervals can be increased to achieve convergence. Though many works have focused on these methods and their application to PDE-constrained OC \cite{JieZhu2023, Heinkenschloss1996, KupferSachs1992, KhaksarShamsi2017, HolmqvistMagnusson2016}, few works have combined them in a unified framework for the general solution of OCPs constrained by PDEs \cite{DaviesPollock2025, DaviesPollock2026}. To the authors' knowledge, \textit{no} work has been done to combine the two methods in an \textit{adaptive} framework for PDE-constrained OC. 

For a method to be adaptive, an estimate of the solution error must be computed, and the mesh must be refined or coarsened in regions that require more or less resolution. A benefit of the use of the finite element method is that there exists a plethora of literature on the derivation of a posteriori error estimates (e.g.~\cite{AinsworthCraig1991,AinsworthOden1997,babuvskaRheinboldt1978,BabuvskaSuri1990, ZienkiewiczZhu1987}), with Ref.~\cite{AinsworthOden1997} classifying a posteriori estimators into three categories: gradient-recovery methods, explicit residual methods, and implicit residual methods.




To provide an estimate for the solution error we employ a residual estimation method, specifically, an implicit residual estimation method, where a series of local residual problems are solved to provide an estimate for the error over an element (or interval). Implicit residual estimators are often avoided as a result of computational expense; however, as the advancement of computational capability has continued to improve, parallelization techniques have allowed for a significant reduction in the computational cost \cite{AinsworthOden1997}. Unlike explicit estimators, implicit estimators circumvent issues arising from dependence on unknown data and sometimes cumbersome derivations of explicit error bounds. Moreover, explicit estimators frequently involve unknown weighting factors associated with boundary and interior errors. As these weights are not easily determined, this limitation further motivates the use of an implicit approach, which eliminates the need for such weighting altogether. 

The implicit residual estimation approach, in addition to its application in the spatial dimension, is extended in a novel procedure to estimate the error associated with the temporal discretization, a so-called $hp-$flipped Legendre-Gauss-Radau ($hp-$fLGR) scheme. The $hp-$fLGR scheme is a popular method for ODE-constrained OCPs known as a \textit{Gaussian quadrature orthogonal collocation} method. Gaussian quadrature orthogonal collocation methods derive the support points within a particular interval from the roots of orthogonal polynomials (e.g.~Legendre, Chebyshev, or Jacobi polynomials). The roots are often referred to as Gaussian quadrature points and are associated with highly-accurate quadrature rules for polynomial functions. For OCPs with smooth solutions, it can be shown that Gaussian quadrature orthogonal collocation methods converge at an exponential rate \cite{HagerHou2019}. The primary benefit of the extension of the implicit residual estimation approach to the $hp-$fLGR method is that the error in both the temporal and spatial dimensions may be equilibrated through the use of comparable error estimators. 
 
Once an estimate of the error on the current mesh has been computed, a new appropriate mesh is generally selected algorithmically through a process known as \textit{mesh refinement}. The solution of the NLP is produced iteratively, where upon each iteration, an estimate of the solution error on the current mesh is computed, and the mesh is locally refined in regions that do not satisfy a requested error tolerance. The decision to employ either $h-$ or $p-$refinement is typically guided by a prediction of the regularity of the solution over each element. In general, it is more advantageous to achieve convergence in nonsmooth regions with $h-$refinement, and in smooth regions with $p-$refinement \cite{DarbyHager2011a}. 

In this paper, the decision on whether to perform $h-$ or $p-$refinement is informed based on the decay rate of the coefficients of a Legendre polynomial expansion of the solution. It has been shown in Refs.~\cite{Mavriplis1994, LiuRao2017, HoustonSenior2003} that the decay rate of the coefficients of a Legendre polynomial approximation of the solution can be used to estimate the regularity of the solution. Larger values of the decay rate indicate smooth profiles and smaller values of the decay rate indicate nonsmooth structure. In addition to a mesh refinement algorithm, a means for mesh reduction is adopted from Ref.~\cite{LiuHager2015} in the temporal and spatial dimension. 

Few direct collocation methods have been developed as adaptive $hp-$methods for the OC of partial differential equation systems. Largely, those that have been developed have focused on their application to the temporal discretization \cite{Betts2020,BettsCampbell2005,BettsCampbell2004, DaviesPollock2026}. No works have been produced to combine adaptive orthogonal collocation and adaptive finite element discretizations in a single, direct framework for accurate solutions to PDE-constrained OCPs. The aim of this work is to extend methods for temporal and spatial adaptivity to a direct collocation approach for the OC of parabolic partial differential equations. 

The novelty and contribution of this work is as follows. A novel $hp-$\textit{adaptive} method for the general solution of OCPs constrained by parabolic partial differential equations is presented. An implicit residual-based error estimator is used to quantify the error in both the temporal and spatial dimensions. The error estimate is novel for orthogonal collocation methods. Based on the provided error estimate, a mesh refinement algorithm is implemented that harnesses the decay rates of Legendre polynomial coefficients to inform $h-$ and $p-$refinement where necessary. This algorithm, previously applied to ODE-constrained OCPs \cite{LiuRao2017} is extended here to PDE-constrained OCPs. Additionally, a means to reduce the mesh size in overcollocated regions that satisfy the mesh tolerance is implemented from Ref.~\cite{LiuHager2015}. Lastly, several numerical examples are performed that highlight the key results and advantages of the $hp-$adaptive finite element method in practice. 

The remainder of this paper is structured as follows. The notations and conventions used in this work are provided in Section \ref{sec:notcov}. An outline of a general OCP governed by a parabolic PDE (and its subsequent discretization) is presented in Section \ref{sec:OCP}. The developed error estimate and error computation procedure is presented in Section \ref{sec:error}. Refinement and reduction procedures based on the calculated error are described in Section \ref{sec:refine}. To demonstrate the effectiveness of the adaptive method, the framework is tested on two numerical examples in Section \ref{sec:examples}. Lastly, conclusions and key findings of the work are summarized in Section \ref{sec:conclusions}. 

\section{Notation and Conventions}\label{sec:notcov}
In this paper, some notational shortcuts are used to aid the reader and maintain conciseness. Relevant notation and conventions are provided here for reference. All vectors are denoted by a bold marking. Matrices are not bolded. Given the vector $\mathbf{f} \in \mathbb{R}^{1 \times n}$ and vector $\mathbf{g} \in \mathbb{R}^{1 \times n}$, an example of a relevant operation would be 
$$
    F = \mathbf{f}^T\mathbf{g} \in \mathbb{R}^{n \times n},
$$
where $\mathbf{f}^T \in \mathbb{R}^{n \times 1}$ denotes the transpose of vector $\mathbf{f}$. Likewise to matrices, scalars are not bolded, a similar operation provides
$$    
f = \mathbf{f}\mathbf{g}^T \in \mathbb{R}. 
$$
To refer to a specific element of a vector, subscript notation is used. That is, if the vector $\mathbf{f} \in \mathbb{R}^{1 \times n}$ is defined by 
$$
    \mathbf{f} = \begin{bmatrix}
        f_1 & f_2 & \ldots & f_n
    \end{bmatrix},
$$
the $i^{\mathrm{th}}$ element of vector $\mathbf{f}$ is denoted by $f_i$. For matrices in this work, a single subscript is used to indicate a specific column of the matrix. Given the matrix $F \in \mathbb{R}^{n \times n}$ written as
$$
    F  = \begin{bmatrix}
        \mathbf{F}_1 & \mathbf{F}_2 & \ldots & \mathbf{F}_n
    \end{bmatrix},
$$
the $i^{\mathrm{th}}$ column of matrix $F$ is denoted by $\mathbf{F}_i \in \mathbb{R}^{n \times 1}$. To refer to a specific subset of indices, a complete list of notation is provided:
\begin{enumerate}
    \item Given the matrix $F \in \mathbb{R}^{n \times n}$, one may extract a specific subset of columns by 
    $$
    F_{i:j} = \begin{bmatrix}
        \mathbf{F}_i & \ldots & \mathbf{F}_j
    \end{bmatrix},
    $$
    where $j \geq i$.
    \item Given the matrix $F \in \mathbb{R}^{n \times n}$, to extract one column of $F$, the result is a vector, thus, we have
    $$
    \mathbf{F}_i \triangleq i^{\mathrm{th}}\:\mathrm{column}\:\mathrm{of}\:F.
    $$
    \item Given the matrix $F \in \mathbb{R}^{n \times n}$, to refer to a specific subset of \textit{rows} of $F$, we first define $F^T$ and write
    \begin{equation*}
        F_{i:j}^T 
        \triangleq i^{\mathrm{th}}\:\mathrm{through}\:j^{\mathrm{th}}\:columns\:\mathrm{of}\:F^T\:(\mathrm{or}\:\mathrm{rows}\:\mathrm{of}\:F).
    \end{equation*}
    \item Given the vector $\mathbf{f} \in \mathbb{R}^{1 \times n}$, one may extract a specific subset of columns by 
    $$
    \mathbf{f}_{i:j}   = \begin{bmatrix}
        f_i & \ldots & f_j
    \end{bmatrix},
    $$
     where $j \geq i$.
     \item Given the vector $\mathbf{f} \in \mathbb{R}^{1 \times n}$, to extract one column of $\mathbf{f}$, the result is a scalar, thus, we have
     $$
     f_i \triangleq i^{\mathrm{th}}\:\mathrm{column}\:\mathrm{of}\:\mathbf{f}.
     $$
\end{enumerate}
For conciseness, derivatives will be expressed with abbreviated notation. For example, the partial derivative of the function $y$ with respect to $x$ can be shortened by 
$$
    \frac{\partial y}{\partial x} \triangleq \partial_xy.
$$
The derivative operator $\partial_x(\cdot)$ only applies to the immediately subsequent variable. If the derivative of an expression of multiple variables is required, parenthesis will be included to indicate the function of differentiation (i.e.~if the derivative of the product of functions $y$ and $g$ is required, we state $\partial_x(yg)$). The total time derivative of the function $y$ can be expressed as 
$$
    \frac{\mathrm{d} y}{\mathrm{d} t} \triangleq \dot{y}.
$$
Further, some important vector and matrix operations are defined as follows. Given the vector $\mathbf{f}(x,t) \in \mathbb{R}^{1 \times n}$, the partial derivative with respect to $x$ is defined as 
$$
    \partial_x \mathbf{f}(x,t) = \begin{bmatrix}
        \partial_x f_1(x,t) & \ldots & \partial_x f_n(x,t)
    \end{bmatrix}.
$$
The integral of $\mathbf{f}(x,t)$ with respect to $x \in [x_0,x_f]$ is defined as 
\begin{equation*}
    \int_{x_0}^{x_f} \mathbf{f}(x,t)\:\mathrm{d}x 
    = \begin{bmatrix}
    \int_{x_0}^{x_f} f_1(x,t)\:\mathrm{d}x &  \ldots & \int_{x_0}^{x_f} f_n(x,t)\:\mathrm{d}x
\end{bmatrix}.
\end{equation*}
If the integral of the matrix $F(x,t) \in \mathbb{R}^{n \times n}$ over $x \in [x_0,x_f]$ is denoted by $G(t)$, the $ij^{\mathrm{th}}$ element of $G$ is defined as 
$$
    G_{ij}(t) = \int_{x_0}^{x_f} F_{ij}(x,t). 
$$
Some vector-valued functions in this paper are functions of other vectors. If $\mathbf{w} \in \mathbb{R}^{1 \times n}$, is a vector of $n$ functions (i.e.~ $\mathbf{w} = [
    w_1(\cdot) \quad w_2(\cdot)\quad \ldots \quad w_n(\cdot)
]$) some possible operations and their resultant outputs are provided below:
\begin{enumerate}
    \item $\mathbf{w}(\mathbf{g},u,t) = \mathbf{r} \in \mathbb{R}^{1 \times n}$, where
    $$ \mathbf{r} = \begin{bmatrix} w_1(g_1,u,t) &  \ldots & w_n(g_n,u,t) \end{bmatrix},$$
    \item $\mathbf{w}(\mathbf{g}^T,u,t) = R \in \mathbb{R}^{n \times n}$, where
    $$ R = \begin{bmatrix} w_1(g_1,u,t) & \ldots & w_n(g_1,u,t)  \\
    \vdots &  \ddots & \vdots  \\
    w_1(g_n,u,t)  & \ldots & w_n(g_n,u,t)  \\
    \end{bmatrix}.$$
\end{enumerate}
All finite-dimensional functions will be presented with dependencies unless specifically stated, whereas vectors of coefficients will be absent of dependencies.

\section{Optimal Control of a Parabolic Partial Differential Equation}\label{sec:OCP}
 Without loss of generality, an OCP is presented in Lagrange form as follows. The goal is to minimize
\begin{equation}
    \mathcal{J} = \int_{\Omega_t}\int_{\Omega_x} \mathcal{L}(x,t,y(x,t))\:\mathrm{d}x\:\mathrm{d}t +\int_{\Omega_t}\mathcal{P}(t,u_1(t),u_2(t))\:\mathrm{d}t, \label{eq:objective}
\end{equation}
subject to the one-dimensional nonlinear PDE  
\begin{equation}
    c_1\partial_t y  +\kappa(y)\partial_x y = c_2 \partial_x(\partial_x y) + f(x,t), 
\end{equation}
with boundary controls, inequality constraints, and fixed-state initial conditions, respectively:
\begin{align}
b_1\left(y(x_0,t),\partial_x y(x_0,t), u_1(t)\right) &= 0, \label{eq:bc1} \\
b_2\left(y(x_f,t),\partial_x y(x_f,t), u_2(t)\right) &= 0, \label{eq:bc2} \\
d\left(y(x,t),x,t,u_1(t),u_2(t)\right) & \leq 0, \label{eq:ineqconst} \\
y(x,t_0)- q(x,t_0) &= 0, \label{eq:initialcondition}
\end{align}
where $\kappa(y)$ is an integrable function of the state variable, $f(x,t)$ is a source function, and $c_1$ and $c_2$ are constant parameters. The problem is defined on $\Omega_t = [t_0, t_f]$ and $\Omega_x = [x_0,x_f]$. A variational form of the PDE can be constructed through the multiplication of the equation by a test function from a suitable space and integration over the spatial domain. The variational form of the PDE is provided by 
\begin{equation}
c_1\int_{\Omega_x} \partial_t yv\:\mathrm{d}x  + \int_{\Omega_x}\kappa(y)\partial_x yv\:\mathrm{d}x = c_2\int_{\Omega_x}\partial_x(\partial_x y)v\:\mathrm{d}x + \int_{\Omega_x}f(x,t)v\:\mathrm{d}x,
\label{eq:variational1}
\end{equation}
where $v$ is the test function defined in $V = H^1(\Omega_x) = \{v\in L_2(\Omega_x) : \partial_x v \in L_2(\Omega_x)\}$ with $L_2(\Omega_x)$ given by the space of square-integrable functions on $\Omega_x$ such that $ \forall g \in L_2(\Omega_x)$:
\begin{equation}
    ||g||_{L_2(\Omega_x)} \triangleq \left( \int_{\Omega_x} |g|^2\:\mathrm{d}x \right)^{\frac{1}{2}} < \infty.
\end{equation}
Integrating Eq.~\eqref{eq:variational1} by parts provides
\begin{multline}
c_1\int_{\Omega_x} \partial_tyv\:\mathrm{d}x  + \int_{\Omega_x}\kappa(y)\partial_xyv\:\mathrm{d}x 
=   \int_{\Omega_x}f(x,t)v\:\mathrm{d}x 
\\+  c_2\left.\partial_xyv\right|_{x = x_f} 
- c_2\left.\partial_xyv\right|_{x = x_0}-c_2\int_{\Omega_x}\partial_xy\partial_xv\:\mathrm{d}x. \label{eq:intbypartsweakform}
\end{multline}
Depending on the type of boundary condition (e.g.~Dirichlet (essential), Neumann (natural), mixed-type), the boundary terms are handled either explicitly through the modification of the $H^1(\Omega_x)$ space or through direct substitution into Eq.~\eqref{eq:intbypartsweakform}. If the boundary conditions are homogeneous Dirichlet conditions (i.e.~$y(x_0,t) = y(x_f,t) = 0$), the subspace $H^1_0(\Omega_x)$ is employed. The space $H^1_0(\Omega_x)$ takes on the same properties of $H^1(\Omega_x)$ with the additional constraint that functions belonging to $H^1_0(\Omega_x)$ vanish in the sense of the trace on $\partial \Omega_x$ (the boundary of $\Omega_x$). When inhomogeneous Dirichlet conditions are present (e.g.~$y(x_0,t) = u_1(t),\:y(x_f,t) = u_2(t)$) one can enforce additional, explicit constraints on the boundary basis functions to assume the value of the Dirichlet data. In fact, this procedure can be identically performed for homogeneous conditions as well; however, the use of the $H^1_0(\Omega_x)$ space implicitly enforces this constraint and thus reduces the required action from the user. If Neumann conditions are assumed, one has
\begin{align}
    \partial_x y(x_0,t)&= g_1(u_1(t)), \label{eq:Neumann1} \\
    \partial_x y(x_f,t) &= g_2(u_2(t)), \label{eq:Neumann2}
\end{align}
where $g_1(u_1(t))$ and $g_2(u_2(t))$ are some functions of the control, we may rewrite Eq.~\eqref{eq:intbypartsweakform} by 
\begin{multline}
c_1\int_{\Omega_x} \partial_tyv\:\mathrm{d}x + \int_{\Omega_x}\kappa(y)\partial_xyv\:\mathrm{d}x 
=   \int_{\Omega_x}f(x,t)v\:\mathrm{d}x
\\ +c_2\left.g_2(u_2(t))v\right|_{x = x_f} 
 - c_2\left.g_1(u_1(t))v\right|_{x = x_0}-c_2\int_{\Omega_x}\partial_xy\partial_xv\:\mathrm{d}x. \label{eq:intbypartsweakformrewrite}
\end{multline}
Thus, the new problem is to find a solution $y$ that satisfies Eq.~\eqref{eq:intbypartsweakformrewrite} for all test functions $v$ in the test space $V$. In this work, the Galerkin approximation is employed, where the solution space (i.e.~the space in which the solution exists) is assumed to be identical to the test space, $V$. The first step under the Galerkin finite element method is to approximate the infinite-dimensional form of the PDE in Eq.~\eqref{eq:intbypartsweakformrewrite} through finite-dimensional functions, $y_h$ and $v_h$, that are constructed by linear combinations of basis functions that define a finite-dimensional subspace of $H^1(\Omega_x)$. These approximations can be written by 
\begin{align}
    y &\approx y_h  = \sum_{k = 1}^K\sum_{j = 1}^{p^{(k)}+1} \phi^{(k)}_j Y^{(k)}_j, \label{findim1} \\
    v &\approx v_h   = \sum_{k = 1}^K\sum_{i = 1}^{p^{(k)}+1} \phi^{(k)}_i V^{(k)}_i, \label{findim2}
\end{align}
where $K$ is the number of elements on the domain $\Omega_x$ and $p^{(k)}$ is the degree, or order, of the $k^\mathrm{th}$ element. In this work P-$p^{(k)}$ Lagrange elements are used, where the governing basis functions within each element are Lagrange polynomials supported at equidistant support nodes. That is, in a particular element, the state approximation may be described by 
\begin{equation}
    y^{(k)} \approx y_h^{(k)} = \sum_{i=1}^{p^{(k)}+1} \phi_j^{(k)} Y_j^{(k)}, \label{eq:stateapproxelement}
\end{equation}
where $\phi_j^{(k)}$ are Lagrange polynomials defined on $r \in [0,1]$ such that 
\begin{equation}
    \phi_j^{(k)} = \begin{cases}
        \prod_{\substack{i=1 \\ j \neq i}}^{p^{(k)}+1} \dfrac{r - r_{i}^{(k)}}{r_{j}^{(k)}-r_{i}^{(k)}}, & x\in \Omega_k, \\
        0, & x \notin \Omega_k,
    \end{cases}
\end{equation}
and $r_i^{(k)}$, $i \in \{1,\ldots,p^{(k)}+1\}$ are a set of strictly monotonically increasing, equidistant support points such that $0 = r_1^{(k)} < r_2^{(k)} < \ldots < r^{(k)}_{p^{(k)}+1} = 1$. The transformed domain is related to the element domain $\Omega_k \in [x_1^{(k)}, x_{p^{(k)}+1}^{(k)}]$ through the associated affine mapping 
\begin{equation}
    x = \left({x^{(k)}_{p^{(k)}+1}-x_1^{(k)}}\right)r + x_1^{(k)}.
\end{equation}
It is advantageous to evaluate the Lagrange polynomials on $r \in [0,1]$ as it provides a consistent domain on which to evaluate the state approximation and its derivatives. It is additionally important to note that the derivative of the Lagrange polynomials on $r$ may be related to the derivative on the physical domain by 
\begin{equation}
   \partial_x \phi_i^{(k)} = \left(\partial_x r\right)^{(k)} \partial_r  \phi_i^{(k)} = \dfrac{1}{h^{(k)}}\partial_r  \phi_i^{(k)}. \label{eq:derivativemapping}
\end{equation}
The Lagrange polynomials $\phi_{p^{(k-1)}+1}^{k-1} $ and $\phi_{1}^{k}$ share the same support point for $k \in\{2,\ldots,K\}$, which is likewise true for $\phi_{p^{(k)}+1}^{k} $ and $\phi_{1}^{k+1}$ for $k \in \{1,\ldots,K-1\}$. To enforce the piecewise continuity constraint at the element boundaries, the additional constraints $Y_{p^{(k-1)}+1}^{(k-1)} = Y_1^{(k)}$ for $k \in \{2,\ldots,K\}$ and $Y_{p^{(k)}+1}^{(k)} = Y_1^{(k+1)}$ for $k \in \{1,\ldots,K-1\}$ are required. However, we may implicitly enforce these constraints by using only one variable in the NLP to represent both coefficients at the boundary of the element. For example, the coefficients $Y_{p^{(k-1)}+1}^{(k-1)} = Y_1^{(k)}$ for a given $k$ may be defined as the same variable in the NLP. This reduces the size of the NLP while implicitly enforcing the required piecewise continuity condition. Thus, the total number of support points (and coefficients) over the entire domain may be computed by 
\begin{equation}
    N_x = \left(\sum_{k = 1}^K p^{(k)}\right) + 1.
\end{equation}
When referring to the state approximation within a given element, the superscript notation will be used to distinguish Eq.~\eqref{eq:stateapproxelement} from Eq.~\eqref{findim1}. Eventually, the finite-dimensional approximation of the PDE is transcribed into a discrete form suitable for numerical implementation. 
\subsection{Nonlinearities in the Optimization Framework}\label{sec:nonlinear}
Nonlinear terms in the variational form are notoriously difficult to address, particularly under a direct collocation framework. A method for transcribing Eq.~\eqref{eq:intbypartsweakformrewrite} into a discrete form that appears linear is provided in Ref.~\cite{DaviesPollock2025}. For brevity, we choose to avoid an in-depth discussion here, though we repeat important details. For more information, the reader is referred to Ref.~\cite{DaviesPollock2025}.


If a new function, $\beta(t)$, is defined by 
\begin{equation}
    \beta(t) = \int_0^t \kappa(s)\:\mathrm{d}s, \label{eq:beta}
\end{equation}
a finite-dimensional weak form can be constructed, written as
\begin{multline}
c_1\int_{\Omega_x} \partial_t y_hv_h\:\mathrm{d}x  + \int_{\Omega_x}\partial_x \beta_hv_h\:\mathrm{d}x 
=   \int_{\Omega_x}f(x,t)v_h\:\mathrm{d}x
\\ +c_2 \left.g_2(u_2(t))v_h\right|_{x = x_f}
- c_2 \left.g_1(u_1(t))v_h\right|_{x = x_0}
-c_2 \int_{\Omega_x}\partial_x y_h\partial_x v_h\:\mathrm{d}x, \label{eq:weakform3}
\end{multline}
where $\beta(y) \approx \beta_h$ and $\beta_h$ is defined explicitly by 
\begin{equation}
    \beta(y) \approx \beta_h = \sum_{k = 1}^{K}\sum_{j = 1}^{p^{(k)}+1} \phi^{(k)}_j \beta^{(k)}_j = \sum_{k = 1}^{K}\sum_{j = 1}^{p^{(k)}+1} \phi^{(k)}_j\int_0^{Y^{(k)}_j} \kappa (s) \:\mathrm{d}s, \label{eq:betaapprox}
\end{equation}
where $\kappa(Y^{(k)}_j)$ is defined by $\kappa(Y^{(k)}_j) \triangleq \left.\kappa(y)\right|_{y = Y^{(k)}_j}$. 
Ultimately, this leads to an approximation of the nonlinear term. The transformation serves as an approximation to the original PDE by neglecting higher-order cross-term products of the basis functions, specifically in the nonlinearity. The ``Kirchoff-like'' transform is a way to improve the computational efficiency in generating a solution to the OCP through a linearized approximation, which can be viewed as an analogy to techniques such as mass lumping \cite{ZienkiewiczTaylor2005, Hughes2003, HintonRock1976} in finite element analysis. One benefit of residual-based error estimation (the focus of this work) is that one may quantifiably demonstrate that the solution of the approximate system approaches the solution of the original nonlinear system as the mesh is refined. The weak form given by Eq.~\eqref{eq:weakform3} has the advantage of being linear in the unknowns $\beta^{(k)}_j$ and $Y^{(k)}_j$, in contrast to the nonlinear dependence on $y$ in Eq.~\eqref{eq:intbypartsweakformrewrite}, which mitigates the computational complexity of evaluating and taking derivatives of the NLP constraints function. If the leading term in the nonlinearity is nonintegrable, other techniques may be required to handle the nonlinearity present. 
\subsection{Construction of the Semi-Discrete Form}
\begin{figure*}[t!]
      \centering
      \includegraphics[scale=1.15]{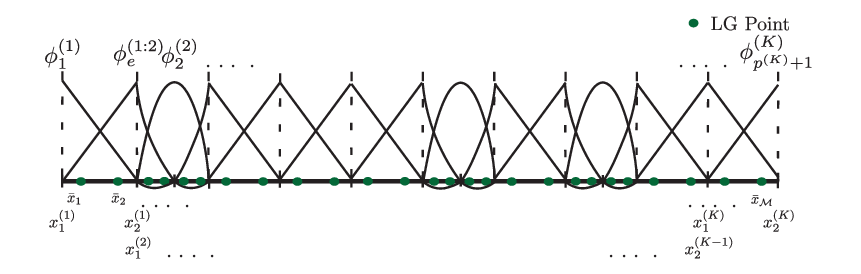}
      \vspace{-1em}
      \caption{A depiction of the spatial discretization.}
      \label{fig:elements}
   \end{figure*}

The semi-discrete form may be constructed by leveraging the linearity of the unknowns $Y^{(k)}_j$ and $\beta^{(k)}_j$. It can be shown that the solution to Eq.~\eqref{eq:weakform3} is independent of the coefficients $V^{(k)}_i$, which yields the transformation of Eq.~\eqref{eq:weakform3} into the matrix-vector semi-discrete form  
\begin{equation}
            c_1M\dot{\mathbf{Y}}(t)+  N\boldsymbol{\beta}(t) = \mathbf{F}(t)  
            -c_2 A\mathbf{Y}(t) 
        + c_2(\mathbf{e}_{N_x}g_2(u_2(t)) - \mathbf{e}_{1}g_1(u_1(t))), \label{eq:semidiscrete}
\end{equation}
where $\boldsymbol{\beta}(t), \mathbf{Y}(t)\in \mathbb{R}^{N_x \times 1}$ are defined by 
\begin{align}
    \boldsymbol{\beta}(t) &= \begin{bmatrix}
    \boldsymbol{\beta}_{1:p^{(1)}}^{(1)} & \boldsymbol{\beta}_{1:p^{(2)}}^{(2)} & \ldots & \boldsymbol{\beta}^{(K)}
    \end{bmatrix}^T, \\
        \mathbf{Y}(t) &= \begin{bmatrix}
    \mathbf{Y}_{1:p^{(1)}}^{(1)} & \mathbf{Y}_{1:p^{(2)}}^{(2)} & \ldots & \mathbf{Y}^{(K)}
    \end{bmatrix}^T,
\end{align}
with $\boldsymbol{\beta}^{(k)},\:\mathbf{Y}^{(k)}\in \mathbb{R}^{1 \times (p^{(k)}+1)}$. If we carefully define the element endpoint basis functions by 
\begin{align}
    \phi_{e}^{(k-1:k)} = \begin{cases}
        \phi^{(k-1)}_{p^{(k-1)}+1}, & x \in \Omega_{k-1}, \\
        \phi^{(k)}_{1}, & x \in \Omega_{k},
    \end{cases} 
\end{align}
$\forall\:k \in \{2,\ldots,K\},$ the vector $\boldsymbol{\phi} \in \mathbb{R}^{1 \times N_x}$ can be constructed as 
\begin{equation}
        \boldsymbol{\phi}\left(r(x,x_1^{(k)},x_{p^{(k)}+1}^{(k)})\right) 
        = \begin{bmatrix}
        \phi_{1}^{(1)} \vspace{1mm} \\ \phi_{\mathcal{I}^{(1)}}^{(1)} \vspace{1mm}\\ \phi_{e}^{(1:2)} \\ \phi_{\mathcal{I}^{(2)}}^{(2)} \vspace{1mm}\\  \phi_{e}^{(2:3)} \\ \vdots \vspace{1mm}\\ \phi_{e}^{(K-1:K)} \\ \phi_{\mathcal{I}^{(K)}}^{(K)} \vspace{1mm}\\ \phi_{p^{(K)}+1}^{(K)} 
    \end{bmatrix}^T,  \label{eq:phivect}
\end{equation}
where $\phi^{(k)}_{\mathcal{I}^{(k)}}$ represents the set of polynomial basis functions associated with the interior of an element
\begin{equation}
    \mathcal{I}^{(k)} \triangleq \{i \in \mathbb{N}\: |\: 2 \leq i \leq p^{(k)}\} \quad k \in \{1,\ldots,K\},
\end{equation}
which allows the matrices $M$, $N$, and $A$ to be computed by  
\begin{align}
    M &= \int_{\Omega_x} \boldsymbol{\phi}^T\boldsymbol{\phi}\:\mathrm{d}x \in \mathbb{R}^{N_x \times N_x}, \label{eq:M} \\
    N &= \int_{\Omega_x} \boldsymbol{\phi}^T\partial_x \boldsymbol{\phi} \:\mathrm{d}x \in \mathbb{R}^{N_x \times N_x}, \label{eq:N} \\
    A &= \int_{\Omega_x} \partial_x \boldsymbol{\phi}^T \partial_x \boldsymbol{\phi}\:\mathrm{d}x \in \mathbb{R}^{N_x \times N_x}. \label{eq:A}
\end{align}
Note that if element $k$ is linear, $\phi_{\mathcal{I}^{(k)}}^{(k)}$ is empty and does not appear in Eq.~\eqref{eq:phivect}. The vector $\mathbf{F}(t) \in \mathbb{R}^{N_x}$ can be computed by 
\begin{equation}
    \mathbf{F}(t) = \int_{\Omega_x} \boldsymbol{\phi}^Tf(x,t)\:\mathrm{d}x,
\end{equation}
and $\mathbf{e}_1$ and $\mathbf{e}_{N_x}$ represent the elementary vectors defined by 
\begin{align}
\mathbf{e}_{1} & = [1,\ldots,0]^T \in \mathbb{R}^{N_x \times 1}, \\
\mathbf{e}_{N_x} & = [0,\ldots,1]^T \in \mathbb{R}^{N_x\times 1}.
\end{align}
It is noted that the integrals in Eqs.~\eqref{eq:M}-\eqref{eq:A} can be computed analytically in certain cases. However, an approximation of the integrals by numerical quadrature (e.g.~Gaussian quadrature) avoids this analytical derivation and generalizes to basis functions of higher-degree. In this work, a Legendre-Gauss (LG) quadrature is utilized. LG quadrature is exact for polynomials of degree $2n-1$, where $n$ is the number of LG points used to approximate the integral. We define $\mathcal{M}$ as the total number of LG points utilized to approximate the spatial integral which can be computed by
\begin{equation}
    \mathcal{M} =\sum_{k = 1}^K 2p^{(k)}.
\end{equation}
It is noted that the number of LG points in each element is twice the degree of the basis functions in the element, ensuring the exactness condition is inherently satisfied. A vector $\mathbf{w} \in \mathbb{R}^{1 \times \mathcal{M}}$ is the set of LG weights on $\Omega_x$. Explicitly, the first $2p^{(1)}$ elements of $\mathbf{w}$ correspond to the weights in element 1, the $2p^{(1)}+1$ through $2p^{(1)}+2p^{(2)}$ elements of $\mathbf{w}$ correspond to the weights in element 2, etc. A spacing vector, $\mathbf{h}\in \mathbb{R}^{1 \times \mathcal{M}}$, is defined similarly where $\mathbf{h}$ is a stacked vector of finite element widths located at each $\bar{x}_p$. For example, the first $2p^{(1)}$ elements of $\mathbf{h}$ are equivalently $x^{(1)}_{p^{(1)}+1} - x^{(1)}_1$ to represent the finite element width at $\bar{x}_{1},\:\bar{x}_2,\ldots,\:\bar{x}_{2p^{(1)}}$.
If we define the matrices $\phi$, $\phi_x \in \mathbb{R}^{N_x \times \mathcal{M}}$ by
\begin{align}
    \phi &= \begin{bmatrix}
        \boldsymbol{\phi}^T [r(\bar{x}_1)] & \ldots & \boldsymbol{\phi}^T [r(\bar{x}_{\mathcal{M}})] 
    \end{bmatrix}, \\
        \phi_x &= \mathbf{h}^{\odot-1} \odot \begin{bmatrix}
        \partial_r \boldsymbol{\phi}^T[r(\bar{x}_1)] &  \ldots & \partial_r \boldsymbol{\phi}^T [r(\bar{x}_{\mathcal{M}})] 
    \end{bmatrix},
\end{align}
where the operator $\mathbf{f} \odot \mathbf{g}$ indicates the Hadamard (element-by-element) product of vectors $\mathbf{f}$ and $\mathbf{g}$, and $\mathbf{h}^{\odot-1}$ indicates the Hadamard inverse (element-by-element reciprocal) of $\mathbf{h}$, we may compute the matrices $M$, $N$, and $A$ numerically by 
\begin{align}
    M &= \left(\frac{1}{2} \odot \mathbf{w} \odot \mathbf{h} \odot \phi \right) \phi^T \in \mathbb{R}^{N_x \times N_x}, \label{eq:Mquad} \\
    N &= \left(\frac{1}{2} \odot\mathbf{w} \odot \mathbf{h} \odot \phi \right) \phi_x^T\in \mathbb{R}^{N_x \times N_x},  \label{eq:Nquad} \\
    A &= \left(\frac{1}{2} \odot\mathbf{w} \odot \mathbf{h} \odot \phi_x \right) \phi_x^T\in \mathbb{R}^{N_x \times N_x}. \label{eq:Aquad} 
\end{align}
 A depiction of the spatial discretization is provided in Fig.~\ref{fig:elements}. The PDE is now removed of explicit spatial dependence, and an application of a temporal discretization is required to fully discretize the problem. 

\subsection{Temporal Discretization}\label{sec:fLGR}
To discretize in time, we employ an orthogonal collocation approach known as the flipped Legendre-Gauss-Radau (fLGR) method. The fLGR method, though not one-to-one with the FEM, displays many similarities in its approach to the discretization. 
The state is supported at a set of fLGR points and a noncollocated initial point. This procedure eliminates the appearance of the Runge phenomenon and provides a useful set of quadrature points for numerical integration. To be specific, fLGR quadrature exhibits the useful property that a polynomial, $p(\tau)$, of degree $2n_t -2$ may be exactly integrated using $n_t$ quadrature points, or more explicitly,
\begin{equation}
    \int_{-1}^{1}p(\tau)\:\mathrm{d}\tau = \sum_{i = 1}^{n_t} w_ip(\tau_i),
\end{equation}
where $w_i$ is the $i^\mathrm{th}$ fLGR quadrature weight. The first step in the fLGR method is to transform the domain. To avoid solving the problem on a dynamic mesh, the time domain is transformed onto a static mesh on $\tau \in [-1,+1]$ by the affine transformation 
\begin{equation}
    \tau = 2\frac{t - t_0}{t_f-t_0} -1, \label{eq:mapping1}
\end{equation}
with 
\begin{equation}
    \dfrac{\mathrm{d} t}{\mathrm{d} \tau} = \dfrac{t_f-t_0}{2}.
\end{equation}
Then, the domain on $\tau \in [-1,+1]$ is broken into a series of $J$ time intervals with $J+1$ monotonically increasing mesh points by $-1 = \tau_1 < \tau_2 <\ldots <\tau_{J+1} = 1$. 
Each mesh interval is transformed onto the domain on which Gaussian quadrature rules exist, $s \in[-1,+1]$, by
\begin{equation}
    s = 2\dfrac{\tau - \tau_{j}}{\tau_{j+1} - \tau_{j}}-1,\quad \{ \tau\:|\: \tau \in [\tau_j,\tau_{j+1}]\}, \label{eq:mapping2}
\end{equation}
$\forall\:j \in \{1,\ldots J\},$ with 
\begin{equation}
    \left(\dfrac{\mathrm{d} \tau}{\mathrm{d}s}\right)^{(j)} = \dfrac{\tau_{j+1} - \tau_{j}}{2}.
\end{equation}
With the appropriate mappings, in the $j^\mathrm{th}$ time interval, we may rewrite $\dot{\mathbf{Y}}(t)$ by 
\begin{equation}
    \dot{\mathbf{Y}}(t) = \dfrac{\mathrm{d} \tau}{\mathrm{d} t}\left(\dfrac{\mathrm{d}s}{\mathrm{d}\tau} \right)^{(j)} \left(\dfrac{\mathrm{d} \mathbf{Y}(s)}{\mathrm{d}s} \right)^{(j)}. \label{eq:derivequiv}
\end{equation}
If we define $\psi^{(j)}$ such that 
\begin{equation}
    \psi^{(j)} = \dfrac{\mathrm{d} t}{\mathrm{d} \tau}\left(\dfrac{\mathrm{d}\tau}{\mathrm{d}s} \right)^{(j)},
\end{equation}
and drop dependencies, Eq.~\eqref{eq:derivequiv} reduces to 
\begin{equation}
        \dot{\mathbf{Y}} = \dfrac{1}{\psi^{(j)}} \left(\dfrac{\mathrm{d} \mathbf{Y}}{\mathrm{d}s} \right)^{(j)}. \label{eq:simpltransform}
\end{equation}
If we substitute Eq.~\eqref{eq:simpltransform} into Eq.~\eqref{eq:semidiscrete}, there exists on each interval 
\begin{equation}
     c_1M\mathbf{Y}_s^{(j)}+  \psi^{(j)}N\boldsymbol{\beta}^{(j)} = \psi^{(j)} \mathbf{F}^{(j)}  -c_2 \psi^{(j)}A\mathbf{Y}^{(j)} 
        + \psi^{(j)}c_2(\mathbf{e}_{N_x}g_2(u_2^{(j)}) - \mathbf{e}_{1}g_1(u_1^{(j)})),  \label{eq:semidiscretetransformed}
\end{equation}
for $j \in \{1,\ldots,J\}$ where $\mathbf{Y}_s^{(j)} = \left(\mathrm{d} \mathbf{Y}/\mathrm{d}s \right)^{(j)}$ and additional constraints are required to ensure continuity in the state at the meshpoints. 

\begin{figure}
    \centering
    \includegraphics[scale = 0.4]{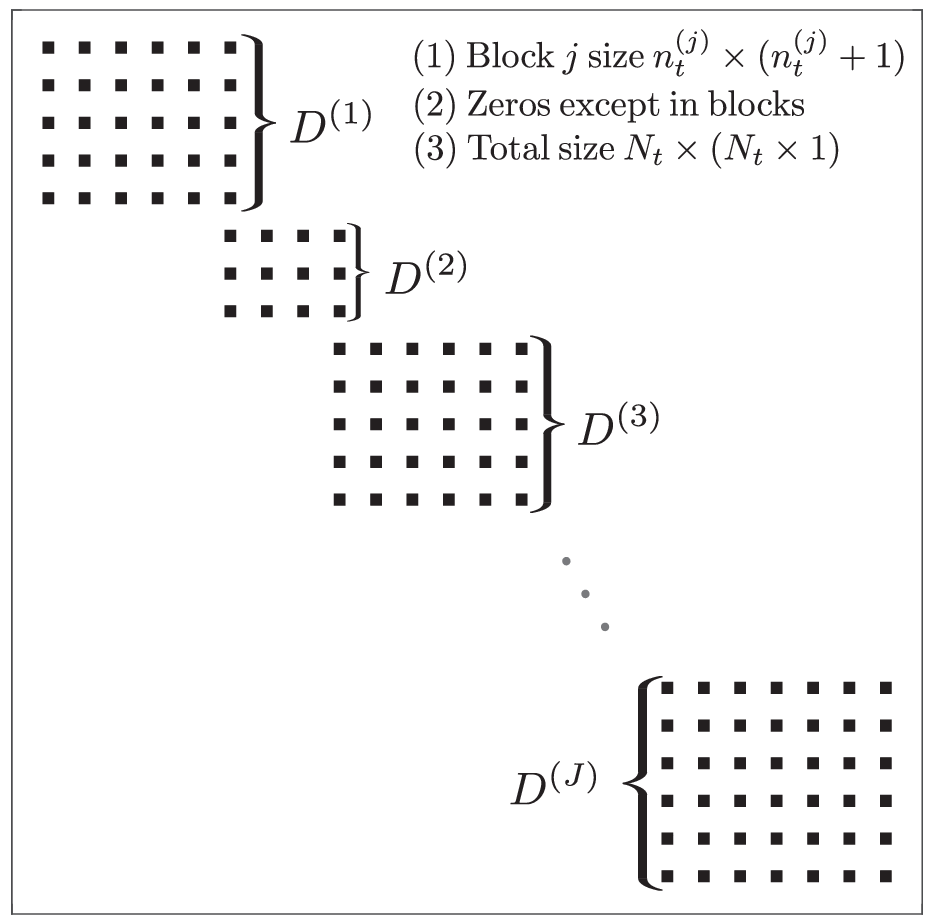}
    \caption{A deptiction of the temporal differentiation matrix, $D_t$, and its organization.}
    \label{fig:diffmatrix}
\end{figure}

The next step in the fLGR method is to approximate the state. The state is approximated by a set of piecewise continuous Lagrange polynomials supported at fLGR points and a noncollocated initial point. The fLGR points are the negative roots of the Legendre polynomial $p_{n_t^{(j)} -1} (s) + p_{n_t^{(j)}}(s)$, where $n_t^{(j)}$ represents the number of fLGR points on interval $j$. The state is approximated in interval $j$ by 
\begin{equation}
    \mathbf{Y}^{(j)} \approx \mathbf{Y}_{h}^{(j)} = \sum_{i = 0}^{n_t^{(j)}} \mathcal{L}^{(j)}_i\mathbf{Y}_i^{(j)}, \label{eq:stateapprox}
\end{equation}
where $\mathbf{Y}_i^{(j)}$ indicates the $i^\mathrm{th}$ column of $Y^{(j)}$, where $Y^{(j)} \in \mathbb{R}^{N_x \times (n_t^{(j)} +1)}$ is defined by 
\begin{equation}
    Y^{(j)} = \begin{bmatrix}
        Y_{10}^{(j)}  & \ldots & Y_{1n_t^{(j)}}^{(j)} \\
        \vdots &  \ddots & \vdots \\
        Y_{N_x0}^{(j)}  & \ldots & Y_{N_xn_t^{(j)}}^{(j)} \\
    \end{bmatrix}. \label{eq:intstatematrix}
\end{equation}
The $i^{\mathrm{th}}$ Lagrange polynomial in interval $j$, $\mathcal{L}^{(j)}_i(s)$, is defined by 
\begin{equation}
      \mathcal{L}^{(j)}_i(s) = \prod_{\substack{k=0 \\ i \neq k}}^{n_t^{(j)}} \dfrac{s - s_{k}^{(j)}}{s_{i}^{(j)}-s_{k}^{(j)}},\label{eq:Lpoly}
\end{equation}
which is of degree $n_t^{(j)}$. The state approximation may be continuously differentiated by 
\begin{equation}
        \dfrac{\mathrm{d} \mathbf{Y}^{(j)}(s)}{\mathrm{d} s} \approx \dfrac{\mathrm{d} \mathbf{Y}_h^{(j)}(s)}{\mathrm{d} s} = \sum_{i = 0}^{n_t^{(j)}} \dfrac{\mathrm{d} \mathcal{L}^{(j)}_i(s)}{\mathrm{d} s} \mathbf{Y}_i^{(j)}, \label{eq:statetimederiv}
\end{equation}
to yield a linear combination of degree $n_t^{(j)} - 1$ polynomials. In order to uniquely define a polynomial of degree $n_t^{(j)} -1$, it must be supported at $n_t^{(j)}$ points. The $n_t^{(j)}$ points chosen to support the derivative of the state approximation are the fLGR points. Note that the noncollocated point exists to uniquely define the state approximation, as it is of one larger degree than the derivative of the state approximation. If the $n_t^{(j)}$ fLGR points are expressed by $s^{(j)}_1,\:s^{(j)}_2,\ldots,\:s^{(j)}_{n_t^{(j)}}$, the derivative of the state approximation evaluated at the fLGR points may be written more succinctly by 
\begin{equation}
    \dfrac{\mathrm{d} \mathbf{Y}_h^{(j)}(s^{(j)}_k)}{\mathrm{d} s} = \sum_{i = 0}^{n_t^{(j)}} \dfrac{\mathrm{d} \mathcal{L}^{(j)}_i(s^{(j)}_k)}{\mathrm{d} s} \mathbf{Y}_i^{(j)} = \sum_{i = 0}^{n_t^{(j)}}D_{ki}^{(j)}\mathbf{Y}_i^{(j)}, 
\end{equation}
$\forall \:k \in \{1,\ldots,n_t^{(j)}\},$ where $D^{(j)}\in \mathbb{R}^{n_t^{(j)} \times \left(n_t^{(j)}+1\right)}$ is the interval differentiation matrix. To enforce continuity in the state variable across meshpoints, the additional constraint must be imposed that 
\begin{equation}
    \mathbf{Y}_{n_t^{(j)}}^{(j)} = \mathbf{Y}_{0}^{(j+1)}, \quad \forall \:j \in \{1,\ldots,J-1\}.
\end{equation}
Although this constraint may be enforced explicitly in the NLP, it is often advantageous to use the same variable in the NLP to implicitly enforce the constraint. Implicitly enforcing the constraint leads to the construction of a temporal differentiation matrix $D_t \in \mathbb{R}^{N_t \times (N_t+1)}$ with overlapping columns for each interval differentiation matrix. The total number of collocation points in the domain is denoted by $N_t$ and is computed by 
\begin{equation}
    N_t = \sum_{j=1}^J n_t^{(j)}.
\end{equation}
A depiction of the temporal differentiation matrix and its organization is found in Fig.~\ref{fig:diffmatrix}.
Alongside the differentiation matrix, state matrices $\bar{Y}$ and $\bar{\beta}$ are organized as
\begin{align}
    \bar{Y} &= \begin{bmatrix}
        Y^{(1)} & Y_{1:n_t^{(2)}}^{(2)} & \ldots & Y_{1:n_t^{(J)}}^{(J)}
    \end{bmatrix} \in \mathbb{R}^{N_x \times \left(N_t+1\right)}, \\
    \bar{\beta} &= \begin{bmatrix}
        \beta_{1:n_t^{(1)}}^{(1)} & \beta_{1:n_t^{(2)}}^{(2)} & \ldots & \beta_{1:n_t^{(J)}}^{(J)}
    \end{bmatrix} \in \mathbb{R}^{N_x \times \left(N_t\right)},
\end{align}
where $\beta^{(j)}_{i:k}$ denotes the $i^\mathrm{th}$ through $k^{\mathrm{th}}$ columns of $\beta^{(j)}$ which is defined similarly to Eq.~\eqref{eq:intstatematrix}. The last variable to explicitly define is the control. The control is defined as a set of discrete variables located at the collocation points. For example, in interval $j$, the control is defined as 
\begin{align}
    u_1^{(j)} &\approx \mathbf{U}_1^{(j)} = \begin{bmatrix}
        U_{1,1}^{(j)}  & \ldots & U_{1,n_t^{(j)}}^{(j)} 
        \end{bmatrix}\in \mathbb{R}^{1 \times n_t^{(j)}}, \\
        u_2^{(j)} &\approx \mathbf{U}_2^{(j)} = \begin{bmatrix}
        U_{2,1}^{(j)} & \ldots & U_{2,n_t^{(j)}}^{(j)}     
    \end{bmatrix}\in \mathbb{R}^{1 \times n_t^{(j)}},
\end{align}
where the notation $U^{(j)}_{i,k}$ indicates the $i^\mathrm{th}$ control variable located at $s^{(j)}_k$. Control vectors $\bar{\mathbf{U}}_1$ and $\bar{\mathbf{U}}_2$ can then be constructed by 
\begin{align}
    \bar{\mathbf{U}}_1 = \begin{bmatrix}
        \mathbf{U}_1^{(1)}  & \ldots & \mathbf{U}_1^{(J)}
    \end{bmatrix} \in \mathbb{R}^{1 \times N_t}, \\
        \bar{\mathbf{U}}_2 = \begin{bmatrix}
        \mathbf{U}_2^{(1)}  & \ldots & \mathbf{U}_2^{(J)}
    \end{bmatrix} \in \mathbb{R}^{1 \times N_t},
\end{align}
which completes the discretization of all variables in the problem. The source term $\mathbf{F}^{(j)}$ may be evaluated at all collocation points by
\begin{equation}
    F^{(j)} = \begin{bmatrix}
        \mathbf{F}^{(j)}(t^{(j)}_1) &  \ldots & \mathbf{F}^{(j)}\left(t^{(j)}_{n_t^{(j)}}\right)
    \end{bmatrix} \in \mathbb{R}^{N_x \times n_t^{(j)}},
\end{equation}
which, for the entire domain, can be combined to yield
\begin{equation}
    F = \begin{bmatrix}
        F^{(1)} & \ldots & F^{(J)}
    \end{bmatrix} \in \mathbb{R}^{N_x \times N_t}.
\end{equation}
It is imperative that the source term be evaluated on the time domain, $t$, rather than the transformed domain, $s$, for accurate computation. For clarity, $t^{(j)}_1 < t^{(j)}_2 < \ldots < t^{(j)}_{n_t^{(j)}}$ represents the location on the \textit{time} domain of $s^{(j)}_1 < s^{(j)}_2 < \ldots < s^{(j)}_{n_t^{(j)}}$, which can be found through the inverse affine mappings of Eqs.~\eqref{eq:mapping1} and \eqref{eq:mapping2}. Now, it is possible to construct a fully-discrete form of Eq.~\eqref{eq:semidiscretetransformed}. The fully-discrete form is given by 
\begin{multline}
     c_1M\left(D_t \bar{Y}^T \right)^T+  \boldsymbol{\psi} \odot (N{\bar{\beta}}) = \boldsymbol{\psi} \odot F  -c_2 \boldsymbol{\psi} \odot (A\bar{Y}_{1:N_t}) 
        \\
        + c_2\boldsymbol{\psi} \odot (\mathbf{e}_{N_x}g_2(\bar{\mathbf{U}}_2) - \mathbf{e}_{1}g_1(\bar{\mathbf{U}}_1)), \label{eq:fullydiscrete}
\end{multline}
where $\boldsymbol{\psi}\in \mathbb{R}^{1 \times N_t}$ is a stacked vector of mappings at the collocation points. Within each mesh interval, the mapping $\psi^{(j)}$ is a constant. So, $\boldsymbol{\psi}$ is constructed such that the proper mapping for each collocation point is applied. For example, in mesh interval 1, the mapping takes on $\psi^{(1)}$, which will appear as the first $n_t^{(1)}$ entries in $\boldsymbol{\psi}$.

\subsection{Discretization of the Objective}
The only remaining procedure is to discretize the objective. We may remove the explicit dependence on the spatial variable by first turning attention to the spatial integral in Eq.~\eqref{eq:objective}. The spatial integral can be discretized by 
\begin{equation}
    \mathbf{I} = \frac{1}{2} \sum_{p = 1}^{\mathcal{M}}w_p h_p\mathcal{L}\left(\bar{x}_p,\mathbf{T}, \left(\boldsymbol{\phi}[r(\bar{x}_p)] \bar{Y}\right)^T\right) \in \mathbb{R}^{(N_t+1) \times 1},
\end{equation}
where if we define
\begin{equation}
    \mathbf{t}^{(j)} = \begin{bmatrix}
        t^{(j)}_0 & \ldots & t^{(j)}_{n_t^{(j)}}
    \end{bmatrix},
\end{equation}
$\mathbf{T} \in \mathbb{R}^{(N_t +1) \times 1}$ can be defined as 
\begin{equation}
    \mathbf{T} = \begin{bmatrix}
        \mathbf{t}^{(1)} & \mathbf{t}^{(2)}_{1:n_t^{(2)}} & \ldots & \mathbf{t}^{(J)}_{1:n_t^{(J)}}
    \end{bmatrix}^T.
\end{equation}
Now, the integral in the temporal dimension can be discretized by
\begin{equation}
    \mathcal{J} \approx \Phi = \sum_{i = 1}^{N_t} \psi_i \omega_i \left(I_i + \mathcal{P}\left(T_i, \bar{U}_{1,i}, \bar{U}_{2,i}\right) \right) \in \mathbb{R},
\end{equation}
where $\boldsymbol{\omega} \in \mathbb{R}^{N_t \times 1} $ is a vector containing the corresponding quadrature weights on $s^{(j)}_1, \ldots, s^{(j)}_{n_t^{(j)}}\:\forall\:j \in \{1,\ldots, J\}$, or, more explicitly:
\begin{equation}
    \boldsymbol{\omega} = \begin{bmatrix}
        \mathbf{w}_t^{(1)} &  \ldots & \mathbf{w}_t^{(J)}
    \end{bmatrix}^T,
\end{equation}
where 
\begin{equation}
    \mathbf{w}_t^{(j)} = \begin{bmatrix}
        w_1^{(j)} & \ldots & w_{n_t^{(j)}}^{(j)}
    \end{bmatrix} \in \mathbb{R}^{1 \times n_t^{(j)}}.
\end{equation}
The quadrature is performed at the collocation points, so it is important to highlight that $I_0$ and $T_0$ are absent from the quadrature. 
\subsection{The Nonlinear Programming Problem}\label{sec:NLPsec}
The OCP has been fully-transcribed into a nonlinear programming problem. The goal is to minimize 
\begin{equation}
\Phi = \sum_{i = 1}^{N_t} \psi_i \omega_i \left(I_i + \mathcal{P}\left(T_i, \bar{U}_{1,i}, \bar{U}_{2,i}\right) \right), \label{eq:NLP1}
\end{equation}
subject to 
\begin{equation}
      c_1M\left(D_t \bar{Y}^T \right)^T+  \boldsymbol{\psi} \odot (N{\bar{\beta}}) = \boldsymbol{\psi} \odot F  -c_2 \boldsymbol{\psi} \odot (A\bar{Y}_{1:N_t}) 
        + c_2\boldsymbol{\psi} \odot (\mathbf{e}_{N_x}g_2(\bar{\mathbf{U}}_2) - \mathbf{e}_{1}g_1(\bar{\mathbf{U}}_1)),
\end{equation}
initial conditions
\begin{equation}
    \mathbf{Y}^{(1)}_0 = q\left(x_{1:N_x},t_0 \right),
\end{equation}
and any state/control inequality constraints 
\begin{equation}
    d\left(\bar{Y} , x_{1:N_x},\mathbf{T},\bar{\mathbf{U}}_1, \bar{\mathbf{U}}_2\right)  \leq 0. \label{eq:NLPend}
\end{equation}
Derivatives of the NLP functions are beyond the focus of this paper; however, more information may be found in Ref.~\cite{DaviesPollock2025}. In this work, the automatic differentiation software $\mathbf{adigator}$ \cite{weinstein2017algorithm} is used to compute derivatives of the NLP objective and constraints functions. Any NLP software may be used to produce a solution to the problem in Eqs.~\eqref{eq:NLP1} - \eqref{eq:NLPend}; however, in this work, the open-source NLP solver \textit{IPOPT} \cite{BieglerZavala2008} is used to produce the solution to the discretized problem. 

\section{Residual-Based Error Estimation}\label{sec:error}
Residual-based error estimation is a concept that stems from finite element theory \cite{AinsworthOden1997,babuvskaRheinboldt1978,BankWeiser1985}. In this paper, it is used in its original context in the framework of the finite element method, but it is also extended to orthogonal collocation for the development of a novel error estimate. The goal of such an approach is to produce a reliable error estimate with similarity in both dimensions. In this way, a refinement process can be implemented such that the mesh error is reduced and equilibrated in both dimensions. 

The residual-based error estimate used in this work is a specific class of residual-based estimation known as the \textit{implicit element (or interval) residual method} \cite{AinsworthOden1997}. In the context of the finite element method, the implicit element residual approach involves the solution of a local initial value problem over each element for a functional representation of the error. Generally, it is assumed that the error exists in a higher-dimensional space than the finite element subspace (i.e.~more informally, that the ``missing'' portion of the solution is attributed to an insufficient representation of the space in which the solution exists). The implicit element residual method has often been practically avoided due to the assumption that the solution of several local finite element problems would be too computationally expensive for an approach to local refinement. However, the advancement of computational hardware and parallelization techniques have mitigated the computational cost of such an approach, and further, the evaluation of the estimate is far more computationally tractable than a global refinement technique. In addition to the computational advancements that facilitate the implicit approach, implicit residual estimation avoids the complex derivation of explicit error estimators that depend upon unknown constant data. In this section, an estimate for the spatial error is derived, and the concept of implicit residual estimation is extended to orthogonal collocation for the derivation of a temporal error estimate. 

\subsection{The Element Residual Method}

\begin{figure}
    \centering
    \includegraphics[scale = 0.4]{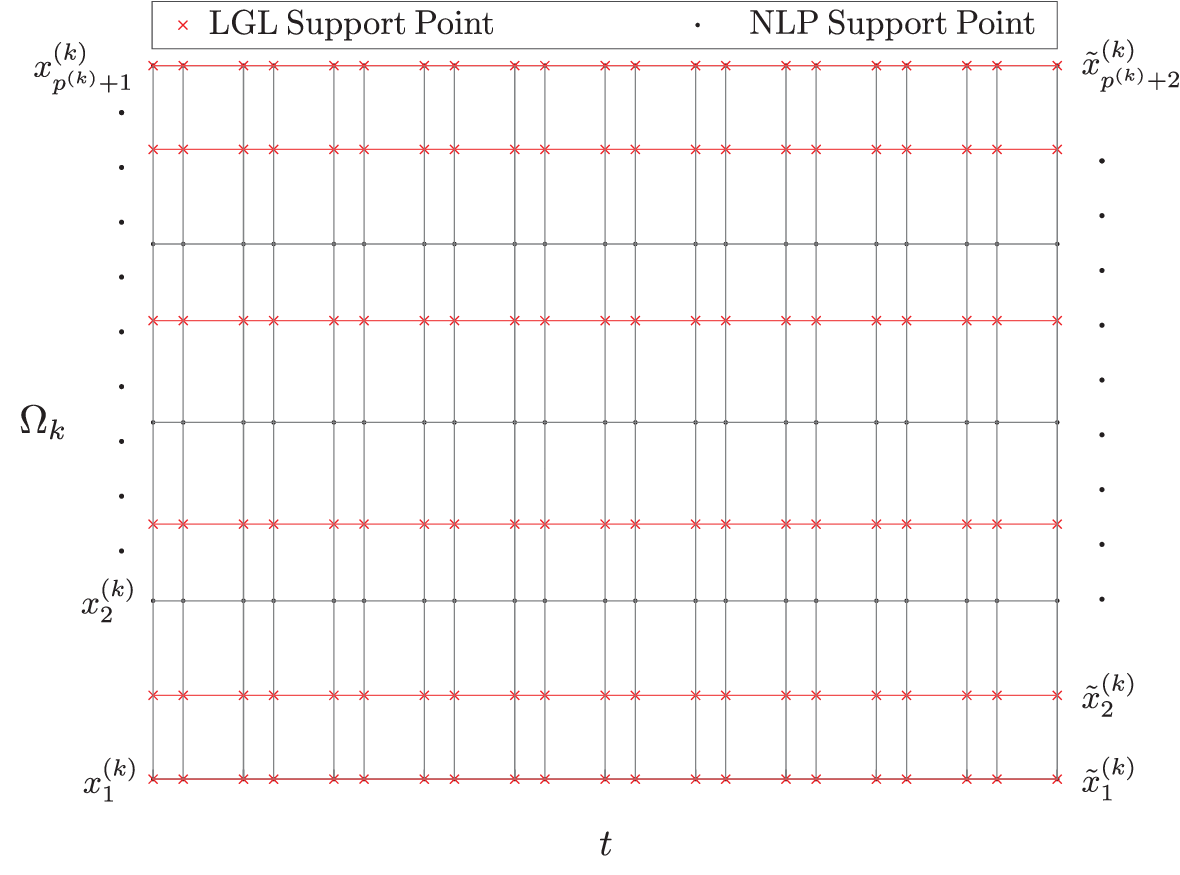}
    \caption{A depiction of the mesh for the local element residual problem with $p^{(k)} = 4$.}
    \label{fig:errormesh}
\end{figure}

The error in the element residual approach is defined as 
\begin{equation}
    e = y -y_h,
\end{equation}
where $y$ is the analytical or ``exact'' solution, which is typically unknown, and $y_h$ is the finite-dimensional finite element solution defined in Eq.~\eqref{findim1}. The first step in the element residual approach is to substitute the definition of the error into the original differential equation:
\begin{equation}
    c_1\partial_t (e+y_h)  +\kappa(e+y_h)\partial_x (e+y_h) 
    = c_2 \partial_x(\partial_x(e+y_h)) + f(x,t).
\end{equation}
For the sake of the remaining discussion, it is assumed that $y_h$ and its derivatives have been computed and are fully determined from the solution to the NLP in Section \ref{sec:NLPsec}. The goal of this approach is to determine a local estimate for the error over an element. As a result, we will multiply the equation by a test function and integrate over the element $\Omega_k \subseteq \Omega_x$ as
\begin{multline}
     \int_{\Omega_k}c_1\partial_t (e^{(k)}+y^{(k)}_h)v^{(k)}\:\mathrm{d}x  
     +\int_{\Omega_k}\kappa(e^{(k)}+y^{(k)}_h) \partial_x (e^{(k)}+y^{(k)}_h)v^{(k)}\:\mathrm{d}x \\ = \int_{\Omega_k}c_2 \partial_x(\partial_x(e^{(k)}+y^{(k)}_h))v^{(k)}\:\mathrm{d}x + \int_{\Omega_k}f(x,t)v^{(k)}\:\mathrm{d}x. \label{eq:errorweakform}
\end{multline}
Now, the problem is to find the unknown error function $e^{(k)}$ in element $k$, $k \in \{1,\ldots,K\}$, for all test functions $v^{(k)} \in V_k = H^1(\Omega_k) = \{v^{(k)}\in L_2(\Omega_k) : \partial_x v^{(k)} \in L_2(\Omega_k)\}$. Again, we require a finite-dimensional approximation of this problem such that it can be managed computationally. It is assumed that the error exists in a higher-dimensional space than the finite element solution $y^{(k)}_h$. That is, if the finite element solution, $y^{(k)}_h$, in element $k$ is of degree $p^{(k)}$, and the solution is supposedly analytic, the error in element $k$ will be attributed to the missing modes of the solution $y^{(k)}$ that are insufficiently captured by a degree $p^{(k)}$ approximation. As a result, we depend upon this assumption through enforcing an approximation of the error function $e^{(k)} \approx e^{(k)}_h$ in element $k$ that is one degree larger than the finite element solution $y^{(k)}_h$. This strategy is common and is found in many works on residual estimation \cite{MitchellMcClain2014,AinsworthOden1997}. Thus, approximations of the form
\begin{align}
    e^{(k)} &\approx e^{(k)}_h = \sum_{j = 1}^{p^{(k)}+2} \Psi^{(k)}_j E^{(k)}_j,  \label{eq:errorapproxspace} \\
    v^{(k)} &\approx v^{(k)}_{h} = \sum_{i = 1}^{p^{(k)}+2} \Psi^{(k)}_i V^{(k)}_i, \label{eq:testapproxspace}
\end{align}
may be constructed where it is noted that as the approximation of the error in element $k$ is of degree $p^{(k)}+1$, each Lagrange polynomial $\Psi^{(k)}$ in the basis must necessarily be supported at $p^{(k)}+2$ points for a unique definition. The support points of the Lagrange polynomial basis must necessarily contain both endpoints of $\Omega_k$, but do not need to maintain the same trend as the original finite element solution (i.e.~equidistant points, Gauss quadrature points, etc.). In this work, a set of $p^{(k)}+2$ Legendre-Gauss-Lobatto (LGL) points are selected as the support points as they provide a convenient set of quadrature points to evaluate the estimate as well as a distinct set of support points from the finite element solution. If $p_n(\xi)$ is the $n^{\mathrm{th}}$ degree Legendre polynomial on $\xi \in [-1,+1]$, the LGL points are the roots of $\dot{p}_{n-1}$ together with the points $-1$ and $+1$. The LGL points are generally associated with the quadrature rule
\begin{equation}
    \int_{-1}^{+1} P(\xi) \:\mathrm{d}\xi \approx \sum_{i =1}^{n} \tilde{w}_i P(\xi_i),
\end{equation}
with LGL weights $\tilde{w}_i$ and LGL points $\xi_i$, $i \in \{1,\ldots, n\}$, where if $P(\xi)$ is a polynomial of degree $2n-3$, the quadrature rule becomes exact. The LGL points on $\xi$ may be mapped onto the original spatial domain through the affine transformation 
\begin{equation}
    \tilde{x}_i^{(k)} = \dfrac{(\xi^{(k)}_i +1)}{2} \left( x^{(k)}_{p^{(k)}+1} - x^{(k)}_{1} \right) + x^{(k)}_{1},
\end{equation}
where $\tilde{x}^{(k)}_i,$ $i \in \{1,\ldots,p+2\}$, refers to the coordinates of the LGL points on the spatial domain. The Lagrange polynomials are defined on $r \in [0,1]$ and can be computed by 
\begin{equation}
    \Psi_j^{(k)}(r(x^{(k)})) = \prod_{\substack{i=1 \\ j \neq i}}^{p^{(k)}+2} \dfrac{r - \tilde{r}_{i}^{(k)}}{\tilde{r}_{j}^{(k)}-\tilde{r}_{i}^{(k)}},
\end{equation}
where $\tilde{r}_{i}^{(k)}$, $i \in \{1,\ldots,p^{(k)}+2\}$, are the LGL points mapped onto $r$. The derivatives of the Lagrange polynomials are evaluated in a similar fashion to the right-hand side of Eq.~\eqref{eq:derivativemapping}; however, for conciseness, the derivatives are expressed in relation to the physical variable rather than the transformed variable. Substituting Eqs.~\eqref{eq:errorapproxspace}, \eqref{eq:testapproxspace} into the weak form in Eq.~\eqref{eq:errorweakform} and integrating by parts provides
\begin{multline}
    c_1\int_{\Omega_k}\partial_t (e^{(k)}_h+y^{(k)}_h)v_h^{(k)}\:\mathrm{d}x 
    +\int_{\Omega_k}\kappa(e^{(k)}_h+y^{(k)}_h) \partial_x (e^{(k)}_h+y^{(k)}_h)v_h^{(k)}\:\mathrm{d}x  = c_2\left.\partial_x (e^{(k)}_h+y^{(k)}_h) v_h^{(k)} \right|_{x = x^{(k)}_{p^{(k)}+1}} \\
    -c_2 \left.\partial_x (e^{(k)}_h+y^{(k)}_h) v_h^{(k)} \right|_{x = x^{(k)}_{1}} 
    - c_2\int_{\Omega_k} \partial_x (e^{(k)}_h+y^{(k)}_h) \partial_x v_h^{(k)}\:\mathrm{d}x 
    + \int_{\Omega_k}f(x,t)v_h^{(k)}\:\mathrm{d}x. \label{eq:errorfindim}
\end{multline}
To complete the discretization of the problem, the boundary conditions and initial conditions must be addressed. The boundary conditions are dependent upon the type of element; that is, whether the element is an interior or boundary element. A boundary element is such that $\Omega_k\:\cap\:\partial \Omega_x \neq 0$; whereas, for an interior element, $\Omega_k\:\cap\:\partial \Omega_x = 0$. On an interior element, the boundary conditions must be approximated. In Eq.~\eqref{eq:errorfindim}, the boundary conditions are conditions on the ``true'' or exact gradient of the solution at the element endpoints, or more explicitly
\begin{align}
    \left.\partial_x (e_h^{(k)} + y_h^{(k)})\right|_{x^{(k)}_{p^{(k)}+1}} &= \left.\partial_x y^{(k)}\right|_{x^{(k)}_{p^{(k)}+1}}, \\ 
    \left.\partial_x (e_h^{(k)} + y_h^{(k)})\right|_{x^{(k)}_{1}} &= \left.\partial_x y^{(k)}\right|_{x^{(k)}_{1}}.
\end{align}
A true gradient at the element endpoints is unknown. At an interior interface, the finite element solution admits distinct left and right derivatives. Since $y_h$ is generally $C^0$ but not $C^1$, these derivatives need not coincide. A numerical interface gradient is therefore defined to combine the two derivatives. At the near interior interface, $x_1^{(k)} = x_{p^{(k-1)}+1}^{(k-1)}$, we can define
\begin{align}
\left.\partial_x y_h\right|_{x_{1}^{(k)-}}
&\triangleq 
\left.\partial_x y_h^{(k-1)}\right|_{x_{p^{(k-1)}+1}^{(k-1)}}, \\
\left.\partial_x y_h\right|_{x_{1}^{(k)+}}
&\triangleq 
\left.\partial_x y_h^{(k)}\right|_{x_{1}^{(k)}},
\end{align}
and at the far interior interface, $x_1^{(k+1)} = x_{p^{(k)}+1}^{(k)}$, we can define
\begin{align}
\left.\partial_x y_h\right|_{x_{p^{(k)}+1}^{(k)-}}
&\triangleq 
\left.\partial_x y_h^{(k)}\right|_{x_{p^{(k)}+1}^{(k)}}, \\
\left.\partial_x y_h\right|_{x_{p^{(k)}+1}^{(k)+}}
&\triangleq 
\left.\partial_x y_h^{(k+1)}\right|_{x_{1}^{(k+1)}}.
\end{align}
If we denote a numerical interface gradient by $\partial y^{(k)}(x_{p^{(k)}+1}^{(k)},\mathbf{T}_c)/\partial x \approx \mathbf{f}_{f}\in \mathbb{R}^{N_t \times 1}$ and $\partial y^{(k)}(x_{1}^{(k)},\mathbf{T}_c)/\partial x \approx \mathbf{f}_{n}\in \mathbb{R}^{N_t \times 1}$, each can be defined, respectively, as  
\begin{align}
    \mathbf{f}_{f} &= \begin{cases}
        \frac{1}{2}
\left(
\left.\partial_x y_h\right|_{x_{p^{(k)}+1}^{(k)-}}
+
\left.\partial_x y_h\right|_{x_{p^{(k)}+1}^{(k)+}}
\right), & k \neq K,  \\
        g(\bar{\mathbf{U}}^T_2), & k = K,
    \end{cases}  \label{eq:fluxfar} \\
      \mathbf{f}_{n} &=  \begin{cases}
         \frac{1}{2}
\left(
\left.\partial_x y_h\right|_{x_{1}^{(k)-}}
+
\left.\partial_x y_h\right|_{x_{1}^{(k)+}}
\right), & k \neq 1, \\
         g(\bar{\mathbf{U}}_1^T), & k = 1,
     \end{cases} \label{eq:fluxnear}
\end{align}
where the vector $\mathbf{T}_c \in \mathbb{R}^{N_t \times 1}$ denotes the portion of $\mathbf{T}$ that exists as a collocation point in the NLP (i.e.~all of $\mathbf{T}$ \textit{except} for $t_0$). Thus, on an interior element, the approximate gradients in Eqs.~\eqref{eq:fluxfar}, \eqref{eq:fluxnear} may be substituted into their respective terms in the weak form in Eq.~\eqref{eq:errorfindim} to form a local Neumann residual problem. On boundary elements, the portion of the boundary of $\Omega_x$ that intersects $\partial \Omega_x$ is handled in a different manner. If the conditions of the original problem are Neumann conditions (in this case, we have assumed them to be), the Neumann data may be substituted explicitly into the appropriate boundary condition in place of the approximation as demonstrated in Eqs.~\eqref{eq:fluxfar}, \eqref{eq:fluxnear}. For the assumed conditions in the NLP in Eqs.~\eqref{eq:NLP1}-\eqref{eq:NLPend}, the conditions would take on the form in Eqs.~\eqref{eq:Neumann1}, \eqref{eq:Neumann2}. If the conditions are Dirichlet conditions, then the basis is modified such that the boundary basis function vanishes at the boundary. That is, 
\begin{align}
    e_h^{(1)}(x_1^{(1)}, \mathbf{T}_c) &= \mathbf{0}, \\
    e_h^{(K)}(x_{p^{(K)}+1}^{(K)}, \mathbf{T}_c) &= \mathbf{0},
\end{align}
as, from the solution to the NLP:
\begin{align}
y^{(1)}(x_1^{(1)},\mathbf{T}_c) &= y_h^{(1)}(x_1^{(1)}, \mathbf{T}_c), \\
    y^{(K)}(x_{p^{(K)}+1}^{(K)}, \mathbf{T}_c) &= y_h^{(K)}(x_{p^{(K)}+1}^{(K)}, \mathbf{T}_c).
\end{align}
For the initial conditions, the error is simply defined by 
\begin{equation}
    e_h^{(k)}(\tilde{x}_i^{(k)},t_0) = q(\tilde{x}_i^{(k)},t_0) - y_h^{(k)}(\tilde{x}_i^{(k)},t_0),  \label{eq:errorics}
\end{equation}
for $i \in \{1,\ldots,p^{(k)}+2\}$. As a result, $K$ element residual problems are constructed, where the problem is to determine the coefficients $E^{(k)}_j$, $j \in \{1,\ldots,p^{(k)}+2\}$, at all points in time. The temporal discretization is held constant from the solution to the original NLP. A depiction of the mesh for an element of degree 4 is provided in Fig.~\ref{fig:errormesh} for reference.
The error function is a measure of the satisfaction (or dissatisfaction) of the dynamical system throughout the element. The benefit of the element residual approach is that each residual problem is a \textit{simulation} problem rather than an \textit{optimization} problem. As a result, derivatives of an NLP constraints function are not required, and the error dynamics may be evaluated on the original nonlinear differential equation with the solution to the NLP at less computational cost. This provides a quantifiable measure for the closeness of the solution of the approximate dynamics in the NLP to the original nonlinear differential equation. 

Each integral term in Eq.~\eqref{eq:errorfindim} may be constructed into a discrete form through a quadrature approximation. This is done in a similar manner to Eqs.~\eqref{eq:M}-\eqref{eq:A}. Again, LG quadrature is used to approximate the integrals in each term in Eq.~\eqref{eq:errorfindim}. In each element, $L = 2(p^{(k)} +1)$, LG points with the associated weights $\hat{w}^{{(k)}}_\ell,\:\ell \in \{1,\ldots,L\}$, are used. The relevant matrices for the linear terms attached to the error function, $M^{(k)},$ $A^{(k)} \in \mathbb{R}^{(p^{(k)}+2)\times (p^{(k)}+2)}$, are computed by 
\begin{align}
    M^{(k)}_{ij} &= \dfrac{h^{(k)}}{2} \sum_{\ell = 1}^{L} \hat{w}^{(k)}_{\ell} \Psi^{(k)}_j(\hat{x}_\ell^{(k)})\Psi^{(k)}_i(\hat{x}_\ell^{(k)}), \\
            A^{(k)}_{ij} &= \dfrac{h^{(k)}}{2} \sum_{\ell = 1}^{L} \hat{w}^{(k)}_{\ell} \dfrac{\partial \Psi^{(k)}_j(\hat{x}_\ell^{(k)})}{\partial x} \dfrac{\partial \Psi^{(k)}_i(\hat{x}_\ell^{(k)})}{\partial x}, 
\end{align}
with additional terms $\mathbf{q}^{(k)}_\ell,\mathbf{T1}^{(k)} _i,\:\mathbf{T2}^{(k)} _i$, $\mathbf{F}^{(k)}_i \in \mathbb{R}^{N_t \times 1}$:
\begin{align}
    \mathbf{q}^{(k)}_\ell &= {\boldsymbol{\psi}^{T\odot-1}}\odot D_ty_h^{(k)}(\hat{x}_\ell^{(k)},\mathbf{T}),  \\
    \mathbf{T1}^{(k)} _i &= \dfrac{h^{(k)}}{2} \sum_{\ell = 1}^L\hat{w}^{(k)}_{\ell} \mathbf{q}^{(k)}_{\ell}\Psi^{(k)}_i(\hat{x}_\ell^{(k)}), \\
    \mathbf{T2}^{(k)} _i &=  \dfrac{h^{(k)}}{2} \sum_{\ell = 1}^L\hat{w}^{(k)}_{\ell} \partial_x y_h^{(k)}(\hat{x}_\ell^{(k)},\mathbf{T}_c)\partial_x \Psi^{(k)}_i(\hat{x}_\ell^{(k)}), \\
    \mathbf{F}^{(k)}_i &=  \dfrac{h^{(k)}}{2} \sum_{\ell = 1}^L\hat{w}^{(k)}_{\ell} f(\hat{x}_\ell^{(k)},\mathbf{T}_c) \Psi^{(k)}_i(\hat{x}_\ell^{(k)}),
\end{align}
for $i \in \{1,\ldots,p^{(k)}+2\}$ where the notation $\mathbf{T1}^{(k)} _i$ is used to denote the $i^{\mathrm{th}}$ column of $T1^{(k)}  \in \mathbb{R}^{N_t \times (p^{(k)}+2)}$, etc., and $\boldsymbol{\psi}^{T\odot-1}$ denotes the Hadamard inverse (element-by-element reciprocal) of $\boldsymbol{\psi}^{T}$. The nonlinear term must be evaluated at each integration point and must be computed on every iteration of the root-finding procedure (due to the dependence upon $e^{(k)} _h$). If we define
\begin{align}
    \bar{\boldsymbol{\kappa}}_\ell &= \kappa\left(e^{(k)}_h(\hat{x}_\ell^{(k)},\mathbf{T}_c) + y^{(k)}_h(\hat{x}_\ell^{(k)},\mathbf{T}_c) \right) \in \mathbb{R}^{N_t \times 1}, \\
    \bar{\mathbf{e}}_\ell &= \partial_x e^{(k)}_h(\hat{x}_\ell^{(k)},\mathbf{T}_c)  + \partial_x y^{(k)}_h(\hat{x}_\ell^{(k)},\mathbf{T}_c) \in \mathbb{R}^{N_t \times 1},
\end{align}
on each root-finding iteration, the nonlinear term, $T3^{(k)} \in \mathbb{R}^{N_t \times (p^{(k)}+2)},$ may be found by
\begin{equation}
    \mathbf{T3}^{(k)} _i =  \dfrac{h^{(k)}}{2} \sum_{\ell = 1}^L\hat{w}^{(k)}_{\ell} \bar{\boldsymbol{\kappa}}_\ell 
    \odot  \bar{\mathbf{e}}_\ell\partial_x \Psi^{(k)}_i(\hat{x}_\ell^{(k)}), 
\end{equation}
where $\mathbf{T3}^{(k)}_i \in \mathbb{R}^{N_t \times 1}$ with $i \in \{1,\ldots,p^{(k)}+2\}.$ It is noted that the linear terms may be computed prior to the root-finding procedure to enhance the computational efficiency of the solution process. Applying the identical temporal discretization to the error function, the error state matrix, $\bar{E}^{(k)} \in \mathbb{R}^{(p^{(k)}+2) \times (N_t +1)}$, may be constructed by 
\begin{equation}
    \bar{E}^{(k)} = \begin{bmatrix}
        E_{10}^{(k)} &  \ldots & E^{(k)}_{1N_t} \\
        \vdots &  \ddots & \vdots \\
        E_{(p^{(k)}+2)0}^{(k)}  & \ldots & E^{(k)}_{(p^{(k)}+2)N_t} \\
    \end{bmatrix}.
\end{equation}
This yields the fully-discrete form of Eq.~\eqref{eq:errorfindim} which can be written as 
\begin{multline}
    c_1 M^{(k)}\left(\boldsymbol{\psi}^{T\odot-1}\odot D_t \bar{E}^{(k)T} \right)^T + c_1T1^{(k)T} + T3^{(k)T}  
    =c_2\mathbf{e}^{(k)}_{p^{(k)}+2} \mathbf{f}_{f}^{(k)T}\\ - c_2\mathbf{e}^{(k)}_{1} \mathbf{f}_{n}^{(k)T}  - c_2A^{(k)}\bar{E}^{(k)}_{1:N_t} -c_2 T2^{(k)T} + F^{(k)T} \in \mathbb{R}^{(p^{(k)}+2) \times N_t}, \label{eq:spaceerrorfullydiscrete}
\end{multline}
with 
\begin{align}
    \mathbf{e}_1^{(k)} &= \begin{bmatrix}
        1 & 0 & \ldots & 0
    \end{bmatrix}^T \in \mathbb{R}^{(p^{(k)}+2) \times 1}, \\
    \mathbf{e}_{p^{(k)}+2}^{(k)} &= \begin{bmatrix}
        0 & 0 & \ldots & 1
    \end{bmatrix}^T \in \mathbb{R}^{(p^{(k)}+2) \times 1}.
\end{align}
There exists $(p^{(k)} +2) \times N_t$ equations in $(p^{(k)} +2) \times (N_t+1)$ unknowns. The additional set of $(p^{(k)} +2)$ constraints stem from the initial conditions, which can be applied by Eq.~\eqref{eq:errorics}. The system in Eqs.~\eqref{eq:errorics}, \eqref{eq:spaceerrorfullydiscrete} can be solved by any root-finding algorithm. 

Once a solution has been computed, we must determine how to use the computed error function. It is noted in Refs.~\cite{AinsworthOden1997, MitchellMcClain2014} that the energy norm of the computed error function can be used to quantify an error indicator in each element. In this work, we choose to employ the $L_2$-norm of the computed error function to describe the error in each element. This is done for multiple reasons. The works in Refs.~\cite{AinsworthOden1997, MitchellMcClain2014} are constrained to elliptic partial differential equations that focus on adaptivity in one or more spatial dimensions. Thus, there is no temporal error to be concerned with and, as a result, equilibration of the error in a temporal and spatial dimension is not relevant to the work. The $L_2$-norm provides an approach that can be applied in both dimensions; whereas, the concept of the energy norm and its extension to the temporal discretization is not apparently obvious. Further, the purpose of the work in this paper is to develop an estimate that may be applied in a general manner. In reality, the error may be quantified in numerous ways once the error function has been computed. Less stringent estimates may be enforced through an absolute value of the integral of the error or an $L_1$-norm, for example, and more stringent estimates may be enforced through an $H^1$-norm or further combinations of $L_p$-norms of the error function and its derivatives to a desired order. Heuristically, the $L_2$-norm strikes a balance between the strictness of the error indicator and the performance of the refinement algorithm. Stricter error estimates can lead to overpolluted mesh sizes, and less stringent estimates can under-resolve the mesh in regions with rapid modes, resulting in poor resolution of important features and potentially unsatisfactory results. As a result, we define the error indicator, $\boldsymbol{\eta}_x^{(k)} \in \mathbb{R}^{(N_t+1) \times 1}$, in element $k$ by 
\begin{equation}
    \boldsymbol{\eta}_x^{(k)} = \left(\int_{\Omega_k}e_h^{(k)}(x,\mathbf{T}) \odot e_h^{(k)}(x,\mathbf{T}) \:\mathrm{d}x \right)^{\odot\dfrac{1}{2}},
\end{equation}
where $\mathbf{g}^{\odot \frac{1}{2}}$ is the Hadamard root (or element-by-element square root) of vector $\mathbf{g}$, and if we define $ \bar{\mathbf{Y}}^{(k)},\:\bar{\mathbf{Y}}_x^{(k)} \in \mathbb{R}^{(N_t+1)  \times 1}$:
\begin{align}
    \bar{\mathbf{Y}}^{(k)} &= \max_{i \in \{1,\ldots,p^{(k)}+2\}} \begin{bmatrix}
        \left| y^{(k)}_h\left(\tilde{x}_i^{(k)}, \mathbf{T}\right) \right|
    \end{bmatrix}, \\
    \bar{\mathbf{Y}}_x^{(k)} &= \max_{i \in \{1,\ldots,p^{(k)}+2\}} \begin{bmatrix}
        \left| \partial_x y^{(k)}_h\left(\tilde{x}_i^{(k)}, \mathbf{T}\right)
   \right| \end{bmatrix},
\end{align}
and the matrix 
\begin{equation}
    \bar{\chi} = \begin{bmatrix}
        \bar{\mathbf{Y}}^{(k)} & \bar{\mathbf{Y}}^{(k)}_x
    \end{bmatrix} \in \mathbb{R}^{(N_t +1) \times 2},
\end{equation}
then a relative parameter vector $\boldsymbol{\chi}_x^{(k)}\in \mathbb{R}^{(N_t+1) \times 1}$ can be constructed 
where
\begin{equation}
    \chi_{x,i}^{(k)} = 1+\max_{j \in \{1,2\}} \begin{bmatrix}
        \bar{\chi}_{ij}
    \end{bmatrix} \quad i \in \{1,\ldots,N_t+1\}. \label{eq:relparam}
\end{equation}
A local, relative error indicator in each element can then be calculated as 
\begin{equation}
    \eta^{(k)}_x = \max \begin{bmatrix}
        \boldsymbol{\eta}^{(k)}_x \odot \boldsymbol{\chi}_x^{(k)\odot-1} \label{eq:relerrorspace}
    \end{bmatrix}.
\end{equation}
After computing $\eta_x^{(k)}$ in each element, the mesh can be locally refined (or reduced) to adjust the error to a desired error tolerance $\epsilon_x$.

\subsection{The Interval Residual Method}
The interval residual approach is based on the ideas of the element residual method for error estimation. This approach is novel in its application to adaptivity for orthogonal collocation methods. Note, this procedure is not limited to its application to partial differential equation systems. Likewise to the element residual approach, the error is defined by 
\begin{equation}
    \mathbf{e} = \mathbf{Y} - \mathbf{Y}_h,
\end{equation}
where $\mathbf{Y} \in \mathbb{R}^{N_x \times 1}$ can be viewed as an exact, or analytical, temporal function that describes the behavior of the spatial coefficients and $\mathbf{Y}_h \in \mathbb{R}^{N_x \times 1}$ is its approximation in Eq.~\eqref{eq:stateapprox}. As we desire a local error estimate, we may reduce the definition to a particular interval $j \in \{1,\ldots,J\}$ by 
\begin{equation}
    \mathbf{e}^{(j)} = \mathbf{Y}^{(j)} - \mathbf{Y}^{(j)}_h.
\end{equation}
The goal of the interval residual approach is to obtain a description of the error in the temporal dimension; as a result, we begin from the semi-discrete form in Eq.~\eqref{eq:semidiscrete} where explicit dependence on the spatial variable has been eliminated. The application of the error definition to the differential equation is similar to the element residual approach. However, due to the ``Kirchoff-like'' integral approximation from Section \ref{sec:nonlinear}, we must take care to handle the application to $\boldsymbol{\beta}^{(j)}(s)$. From Eq.~\eqref{eq:betaapprox}, we have that 
\begin{equation}
    \beta^{(j)}_k = \int_{0}^{Y^{(j)}_k} \kappa^{(j)}(l)\:\mathrm{d}l, \quad \forall \:k \in \{1,\ldots,N_x\},
\end{equation}
where $l$ is a dummy variable for integration. Thus, the error in $\boldsymbol{\beta}^{(j)}(s)$, $\boldsymbol{\beta}^{(j)}_e(s) \in \mathbb{R}^{N_x \times 1}$, can be found by 
\begin{equation}
    \beta^{(j)}_{e,k} = \int_{0}^{e^{(j)}_k+Y^{(j)}_{h,k}} \kappa(l)\:\mathrm{d}l, \quad \forall \:k \in \{1,\ldots,N_x\}.
\end{equation}
Substituting the error definition into the dynamical equation in Eq.~\eqref{eq:semidiscrete} yields 
\begin{multline}
    c_1M\left[\frac{\mathrm{d}\mathbf{e}^{(j)}(s)}{\mathrm{d} s} + \frac{\mathrm{d} \mathbf{Y}^{(j)}_{h}(s)}{\mathrm{d} s} \right] + \psi^{(j)}N\boldsymbol{\beta}^{(j)}_e(s)  = \psi^{(j)}\mathbf{F}^{(j)}(s)\\-c_2\psi^{(j)}A[\mathbf{e}^{(j)}(s) + \mathbf{Y}^{(j)}_h(s)] 
    +c_2 \psi^{(j)}(\mathbf{e}_{N_x}g_2(u^{(j)}_2(s)) - \mathbf{e}_{1}g_1(u^{(j)}_1(s))). \label{eq:errortimedynamics}
\end{multline}
As with the element residual approach, if the solution is sufficiently regular, the error $\mathbf{e}^{(j)}(s)$ reflects the inability of $\mathbf{Y}^{(j)}_h(s)$ to represent the unresolved modes of $\mathbf{Y}^{(j)}(s)$. As a result, we assume that the error, $\mathbf{e}^{(j)}(s)$, lies in a higher-dimensional space than the NLP solution, where it is assumed that the ``dominant'' missing mode (i.e.~the most significant contributor to the error, $\mathbf{e}^{(j)}(s)$) is the $(n_t^{(j)}+1)^{\mathrm{th}}$ degree term. 

Similar assumptions are not uncommon, such as the one found in Ref.~\cite{PattersonHager2015}, which is also employed in Ref.~\cite{LiuRao2017}. In Ref.~\cite{PattersonHager2015}, a regularity assumption is made to claim that an additional LGR point in each interval yields a state solution of higher accuracy than the NLP solution. Our estimate differs from that found in Ref.~\cite{PattersonHager2015}, as here, we are solving for an unknown function that represents the residual in the dynamical equations throughout the interval. In Ref.~\cite{PattersonHager2015}, two numerical solutions are compared, where it is assumed that one is of greater accuracy than the other. We choose to solve $J$ interval residual problems, where in each interval, the problem is solved at $n_t^{(j)} + 1$ collocation points. 

\begin{figure}[h]
    \centering
    \includegraphics[scale = 0.4]{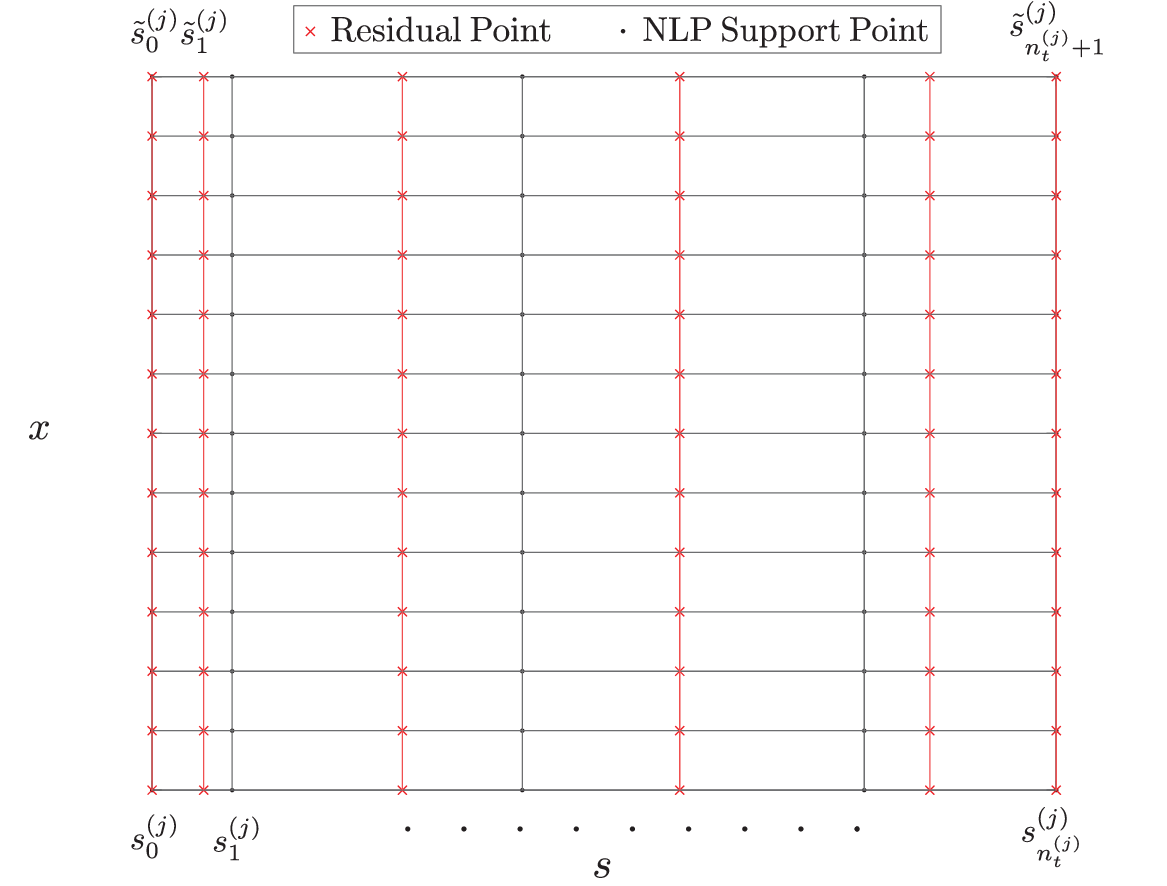}
    \caption{A depiction of the mesh for the interval residual problem.}
    \label{fig:intervalresid}
\end{figure}

The collocation points of the interval residual problem are $n_t^{(j)} + 1$ fLGR points. The $n_t^{(j)} + 1$ collocation points of the interval residual problem form a distinct set of $n_t^{(j)} + 1$ points from the original collocation points of the NLP. Likewise to the temporal discretization of the NLP, a noncollocated support point is included at the start of each interval for consistency with the associated polynomial degree. A depiction of the mesh for the interval residual problem for an interval of degree 4 is provided in Fig.~\ref{fig:intervalresid} for reference. 
The error function in each interval is supported at $n_t^{(j)}+2$ points and is a linear combination of polynomials of degree $n_t^{(j)}+1$. The finite-dimensional form of the error function is given by 
\begin{equation}
    \mathbf{e}^{(j)}(s) \approx \mathbf{e}_h^{(j)}(s) = \sum_{i = 0}^{n_t^{(j)}+1} \ell_i^{(j)}(s)\mathbf{E}^{(j)}_i,
\end{equation}
where $\ell_i^{(j)}(s)$ are the Lagrange polynomials of the interval residual problem defined by 
\begin{equation}
     \ell^{(j)}_i(s) = \prod_{\substack{k=0 \\ i \neq k}}^{n_t^{(j)}+1} \dfrac{s - \tilde{s}_{k}^{(j)}}{\tilde{s}_{i}^{(j)}-\tilde{s}^{(j)}_{k}},
\end{equation}
where $\tilde{s}_i^{(j)}$, $i \in \{0,\ldots,n_t^{(j)}+1\}$, are the $n_t^{(j)}+1$ fLGR points and the noncollocated initial point, $\tilde{s}_0^{(j)}$. The $i^{\mathrm{th}}$ column of $E^{(j)}\in \mathbb{R}^{N_x \times (n_t^{(j)}+2)}$ is denoted by $\mathbf{E}^{(j)}_i$, where $E^{(j)}$ can be written as 
\begin{equation}
  E^{(j)} =   \begin{bmatrix}
    E^{(j)}_{10} &  \ldots & E^{(j)}_{1(n_t^{(j)}+1)} \\
    \vdots & \ddots & \vdots \\
    E^{(j)}_{N_x0} & \ldots & E^{(j)}_{N_x(n_t^{(j)}+1)} \\
\end{bmatrix}.
\end{equation}
The derivative of the error approximation can be computed by 
\begin{equation}
    \frac{\mathrm{d} \mathbf{e}^{(j)}(s)}{\mathrm{d} s} \approx \frac{\mathrm{d} \mathbf{e}_h^{(j)}(s)}{\mathrm{d} s} = \sum_{i = 0}^{n_t^{(j)}+1} \frac{\mathrm{d} \ell_i^{(j)}(s)}{\mathrm{d} s}\mathbf{E}^{(j)}_i,
\end{equation}
where $\hat{D}_{ki}^{(j)}$ can be defined as 
\begin{equation}
    \hat{D}_{ki}^{(j)} = \dfrac{\mathrm{d} \ell_i^{(j)}(\tilde{s}^{{(j)}}_k)}{\mathrm{d} s}.
\end{equation}
To complete the discretization of the interval residual problem, we require an interpolation (or extrapolation) of the control information to the collocated points of the interval residual problem. As only an explicit functional form of the state is present, we require a functional form of the control to interpolate to the collocation points of the residual problem. The interpolant is defined by 
\begin{align}
    U^{(j)}_{1,h} (s) =& \sum_{i=1}^{n_t^{(j)}} \hat{\mathcal{L}}^{(j)}_i(s)U_{1,i}^{(j)}, \label{eq:controlinterp1} \\
    U^{(j)}_{2,h}(s) =& \sum_{i=1}^{n_t^{(j)}} \hat{\mathcal{L}}^{(j)}_i(s)U_{2,i}^{(j)}, \label{eq:controlinterp2}
\end{align}
where $\hat{\mathcal{L}}^{(j)}_i(s)$ is defined by 
\begin{equation}
      \hat{\mathcal{L}}^{(j)}_i(s) = \prod_{\substack{k=1 \\ i \neq k}}^{n_t^{(j)}} \dfrac{s - s_{k}^{(j)}}{s_{i}^{(j)}-s_{k}^{(j)}},\label{eq:LpolyE}
\end{equation}
and it is noted that the initial, noncollocated point of the NLP problem is not included in the interpolant as it does not contain a value for the control. The vectors $\tilde{\mathbf{U}}_1^{(j)} \in \mathbb{R}^{1\times (n_t^{(j)}+1)}$ and $\tilde{\mathbf{U}}_2^{(j)} \in \mathbb{R}^{1\times (n_t^{(j)}+1)}$ represent the evaluation of the finite-dimensional representation of the controls $ U^{(j)}_{1,h} (s)$ and $ U^{(j)}_{2,h} (s)$, respectively, at the collocation points of the residual problem. The definitions can be written explicitly by
\begin{align}
    \tilde{\mathbf{U}}^{(j)}_1 &= \begin{bmatrix}
        U^{(j)}_{1,h} \left(\tilde{s}_1^{(j)}\right) & \ldots & U^{(j)}_{1,h} \left(\tilde{s}_{n_t^{(j)}+1}^{(j)}\right)
    \end{bmatrix}, \\
    \tilde{\mathbf{U}}^{(j)}_2 &= \begin{bmatrix}
        U^{(j)}_{2,h} \left(\tilde{s}_1^{(j)}\right) & \ldots & U^{(j)}_{2,h} \left(\tilde{s}_{n_t^{(j)}+1}^{(j)}\right)
    \end{bmatrix}.
\end{align}
If we define
\begin{equation}
    \tilde{\beta}^{(j)}_{e,ki} = \int_0^{E^{(j)}_{ki} + Y_{ki}^{(j)}} \kappa(l) \:\mathrm{d}l,
\end{equation}
we may complete the discretization of the problem with the matrices $\tilde{Y}^{(j)},\:\tilde{Y}_s^{(j)} ,\:\tilde{F}^{(j)},\:\tilde{\beta}_e^{(j)} \in \mathbb{R}^{N_x \times (n_t^{(j)} +1)}$:
\begin{align}
    \tilde{Y}^{(j)} &= \begin{bmatrix}
        \mathbf{Y}_h^{(j)}\left(\tilde{s}^{(j)}_1\right) &  \ldots &\mathbf{Y}_h^{(j)}\left(\tilde{s}^{(j)}_{n_t^{(j)}+1}\right)
    \end{bmatrix}, \label{eq:Yresidmat}\\
    \tilde{Y}_s^{(j)} &= \begin{bmatrix}
        \dfrac{\mathrm{d} \mathbf{Y}_h^{(j)}\left(\tilde{s}^{(j)}_1\right)}{\mathrm{d} s} & \ldots &\dfrac{\mathrm{d} \mathbf{Y}_h^{(j)}\left(\tilde{s}^{(j)}_{n_t^{(j)}+1}\right)}{\mathrm{d} s}
    \end{bmatrix}, \\
    \tilde{F}^{(j)} &= \begin{bmatrix}
        \mathbf{F}^{(j)}\left(\tilde{s}^{(j)}_1\right) & \ldots &\mathbf{F}^{(j)}\left(\tilde{s}^{(j)}_{n_t^{(j)}+1}\right)
    \end{bmatrix}, \\
    \tilde{\beta}_e^{(j)} &= \begin{bmatrix}
        \tilde{\beta}^{(j)}_{e,11} &  \ldots & \tilde{\beta}^{(j)}_{e,1(n_t^{(j)}+1)} \\
        \vdots  & \ddots & \vdots \\
        \tilde{\beta}^{(j)}_{e,N_x1} &  \ldots & \tilde{\beta}^{(j)}_{e,N_x(n_t^{(j)}+1)}
    \end{bmatrix}. \label{eq:betaresidmat}
\end{align}
The matrices defined in Eqs.~\eqref{eq:Yresidmat}-\eqref{eq:betaresidmat} allow the construction of a fully-discrete form of the interval residual problem in Eq.~\eqref{eq:errortimedynamics}, which can be written as 
\begin{multline}
     c_1M\left[(\hat{D}^{(j)}E^{(j)T})^T+ \tilde{Y}_s^{(j)}\right] + \psi^{(j)}N\tilde{\beta}_e^{(j)} = \psi^{(j)}\tilde{F}^{(j)} \\
     -c_2\psi^{(j)}A[\mathbf{E}^{(j)}_{(1:n_t^{(j)}+1)} + \tilde{Y}^{(j)}]  + c_2\psi^{(j)}\left(\mathbf{e}_{N_x}g_2(\tilde{\mathbf{U}}_2^{(j)}) - \mathbf{e}_{1}g_1(\tilde{\mathbf{U}}_1^{(j)}) \right) \in \mathbb{R}^{N_x \times ( n_t^{(j)} +1 )}. \label{eq:errortimediscrete}
\end{multline}

Equation (\ref{eq:errortimediscrete}) is a system of $N_x \times (n_t^{(j)}+1)$ equations with $N_x \times (n_t^{(j)}+2)$ unknowns, $E^{(j)}$. Imposing the initial condition $\mathbf{E}^{(j)}_0 = \mathbf{0} \in \mathbb{R}^{N_x}$ renders the system completely determined. As the boundary conditions were assumed to be Neumann conditions, the boundary conditions appear explicitly in the dynamical system in Eq.~\eqref{eq:errortimediscrete}. If the conditions are Dirichlet conditions, such as 
\begin{align}
    y(x_0,t) &= b_1(u_1,t), \\
    y(x_f,t) &= b_2(u_2,t),
\end{align} 
the conditions for the interval residual problem can be found by 
\begin{align}
    \mathbf{E}^{(j)T}_{1} &= b_1\left(\tilde{\mathbf{U}}_1^{(j)},s(\tau(t,t_0,t_f),\tau_{j},\tau_{j+1})\right) -  \tilde{\mathbf{Y}}^{(j)T}_{1}, \\ 
    \mathbf{E}^{(j)T}_{N_x} &= b_2\left(\tilde{\mathbf{U}}_2^{(j)},s(\tau(t,t_0,t_f),\tau_{j},\tau_{j+1})\right) -  \tilde{\mathbf{Y}}^{(j)T}_{N_x},
\end{align}
which results in direct constraints on the coefficients of the basis functions on the near and far boundaries and a reduction in the size of Eq.~\eqref{eq:errortimediscrete} by $2 \times ( n_t^{(j)} +1 )$ (i.e.~the equations pertaining to the near and far boundaries are omitted). 

The interval residual problem can now be solved by any root-finding technique such as MATLAB's $\mathbf{fsolve}()$. Once a solution has been computed, we must now decide what to do with the computed error function. To effectively relate the error in the temporal dimension to the error in the spatial dimension, we choose to evaluate the error by an $L_2$-norm over the interval. An error indicator, $\boldsymbol{\eta}_t^{(j)} \in \mathbb{R}^{N_x \times 1}$, can be computed by 
\begin{equation}
    \boldsymbol{\eta}_t^{(j)} = \left(\psi^{(j)}\int_{-1}^{+1} \mathbf{e}_h^{(j)}(s) \odot \mathbf{e}_h^{(j)}(s) \:\mathrm{d}s \right)^{\odot \dfrac{1}{2}},
\end{equation}
and if we define $\hat{\mathbf{Y}}^{(j)},\:\hat{\mathbf{Y}}_t^{(j)} \in \mathbb{R}^{N_x \times 1}$:
\begin{align}
    \hat{\mathbf{Y}}^{(j)} &= \max_{k \in \{0,\ldots,n_t^{(j)}+1\}} \begin{bmatrix}
       \left| \mathbf{Y}^{(j)}_h\left(\tilde{s}_k^{(j)}\right) \right|
    \end{bmatrix}, \\
    \hat{\mathbf{Y}}_t^{(j)} &= \max_{k \in \{0,\ldots,n_t^{(j)}+1\}} \begin{bmatrix}
        \left| \dfrac{1}{\psi^{(j)}} \dfrac{\mathrm{d} \mathbf{Y}^{(j)}_h\left(\tilde{s}_k^{(j)}\right)}{\mathrm{d} s} \right|
    \end{bmatrix},
\end{align}
and the matrix 
\begin{equation}
    \hat{\chi} = \begin{bmatrix}
        \hat{\mathbf{Y}}^{(j)} & \hat{\mathbf{Y}}_t^{(j)}
    \end{bmatrix} \in \mathbb{R}^{N_x \times 2},
\end{equation}
a relative parameter vector $\boldsymbol{\chi}_t^{(j)}\in \mathbb{R}^{N_x \times 1}$ can be constructed 
where 
\begin{equation}
    \chi_{t,i}^{(j)} = 1+\max_{j \in \{1,2\}}\begin{bmatrix}
        \hat{\chi}_{ij}
    \end{bmatrix}\quad i \in \{1,\ldots,N_x\}.
\end{equation}
A local, relative error indicator in each interval can then be calculated as 
\begin{equation}
    \eta^{(j)}_t = \max \begin{bmatrix}
        \boldsymbol{\eta}^{(j)}_t \odot \boldsymbol{\chi}_t^{(j)\odot-1}
    \end{bmatrix}.
\end{equation}
After computing $\eta_t^{(j)}$ in each interval, the mesh can be locally refined (or reduced) to adjust the error to a desired error tolerance $\epsilon_t$. 


\section{Mesh Refinement and Reduction}\label{sec:refine}


It has been shown in Refs.~\cite{Mavriplis1994,LiuRao2017,MitchellMcClain2014}, that the decay rate of the coefficients in a Legendre polynomial expansion can be used to approximate the regularity of the solution. The regularity of the solution can be used to inform the refinement process. The expansion is given in the form 
\begin{equation}
    y(x) \approx y_\ell(x) = \sum_{i = 0}^{\infty} a_i P_i(x),
\end{equation}
where $a_i$ are the expansion coefficients and $P_i(x)$ are the $i^{th}$ degree Legendre polynomials on element $k$. It has been shown in Ref.~\cite{WangXiang2012} that, provided the solution is analytic in a neighborhood of the element, the Legendre polynomial coefficient values $a_i$ decay like $ci10^{-\tilde{\sigma}},\:\tilde{\sigma} >0$. As in Ref.~\cite{LiuRao2017}, we assume a slightly smaller value, $\sigma$, satisfying $ci10^{-\tilde{\sigma}i} \leq c10^{-\sigma i}$ which simplifies the analysis. Thus, the coefficients can be approximated by an exponential least-squares fit of the form
\begin{equation}
    a_i \approx c10^{-\sigma i},\quad \sigma > 0. \label{eq:coeffmodel}
\end{equation}
This, in practice, leads to a more conservative upper bound on the error of the expansion approximation, derived in Ref.~\cite{LiuRao2017} as 
\begin{equation}
\hat{e} = \dfrac{c10^{-\sigma(p^{(k)}+1)}}{\sqrt{1-10^{-2\sigma}}}. \label{eq:errorlegendre}
\end{equation}
By observation, the decay rate of the coefficients $a_i$ as a function of the index, $i$, is identical to the decay rate of the upper bound on the error, $\hat{e},$ as a function of the element degree, $p$. 

In practice, the coefficients, $a_i$, and the corresponding decay rate are unknown and must be estimated from relevant solution data. The Legendre coefficients on element $k$ can be computed by transforming the corresponding finite element solution. From Eq.~\eqref{findim1}, we have that the finite element solution on element $T$ may be written as 
\begin{equation}
    y^{(k)} \approx y_h^{(k)} = \sum_{i=1}^{p^{(k)}+1} \phi^{(k)}_i Y^{(k)}_i,
\end{equation}
where the Lagrange polynomials, $\phi^{(k)}_i$, are supported at $p^{(k)}+1$ equidistant points. It is noted that, as the Legendre polynomials are defined on $\xi \in [-1,+1]$, it is necessary to transform the element's domain such that the support points are appropriately mapped to $\xi$.  We may then write the transformation as  
\begin{equation}
    \begin{bmatrix}
        \hat{a}_0 \\
        \vdots \\
        \hat{a}_p
    \end{bmatrix} = \begin{bmatrix}
        P_0(\xi_1) &  \ldots & P_p(\xi_1) \\
        \vdots & \ddots & \vdots \\
        P_0(\xi_{p+1}) & \ldots & P_p(\xi_{p+1})
    \end{bmatrix}^{-1} \begin{bmatrix}
        Y^{(k)}_1 \\
        \vdots \\
        Y^{(k)}_{p+1}
    \end{bmatrix},
\end{equation}
where the coefficients $Y^{(k)}_i,\:i \in \{1,\ldots,p^{(k)}+1\}$, are the coefficients of the Lagrange polynomial approximation utilized on element $k$, and $\xi_i$, $i \in \{1,\ldots,p^{(k)}+1\}$, are the support points of the Lagrange polynomial approximation mapped onto $[-1,+1]$. The notation $\hat{a}_i$ is used to indicate that the coefficients are computed from the solution data and are approximate.  The coefficients $\hat{a}_i$, $i \in \{0,\ldots,p^{(k)}\}$ are then used in an orthogonal least squares fit to the form in Eq.~\eqref{eq:coeffmodel} to compute the corresponding decay rate $\sigma$. For some functions, for example even or odd functions defined on $[-1,+1]$ \cite{HoustonSenior2003}, it is possible that some of the coefficient values may be zero. In practice, these values must be ``filtered-out'' of the fitting procedure prior to estimating the decay rate. Larger values of the decay rate, say above some $\bar{\sigma}$, indicate that the function is regular, and convergence is best achieved by increasing the degree of the element (i.e.~$p-$refinement). Smaller values (i.e.~below $\bar{\sigma}$) indicate a nonsmooth solution, and convergence is best achieved by partitioning the element into multiple elements (i.e.~$h$-refinement). 

\begin{figure*}[h!]
    \subfloat[$\mathrm{exp}(\xi)$ \label{fig:coeffs1}]{\includegraphics[scale=0.4]{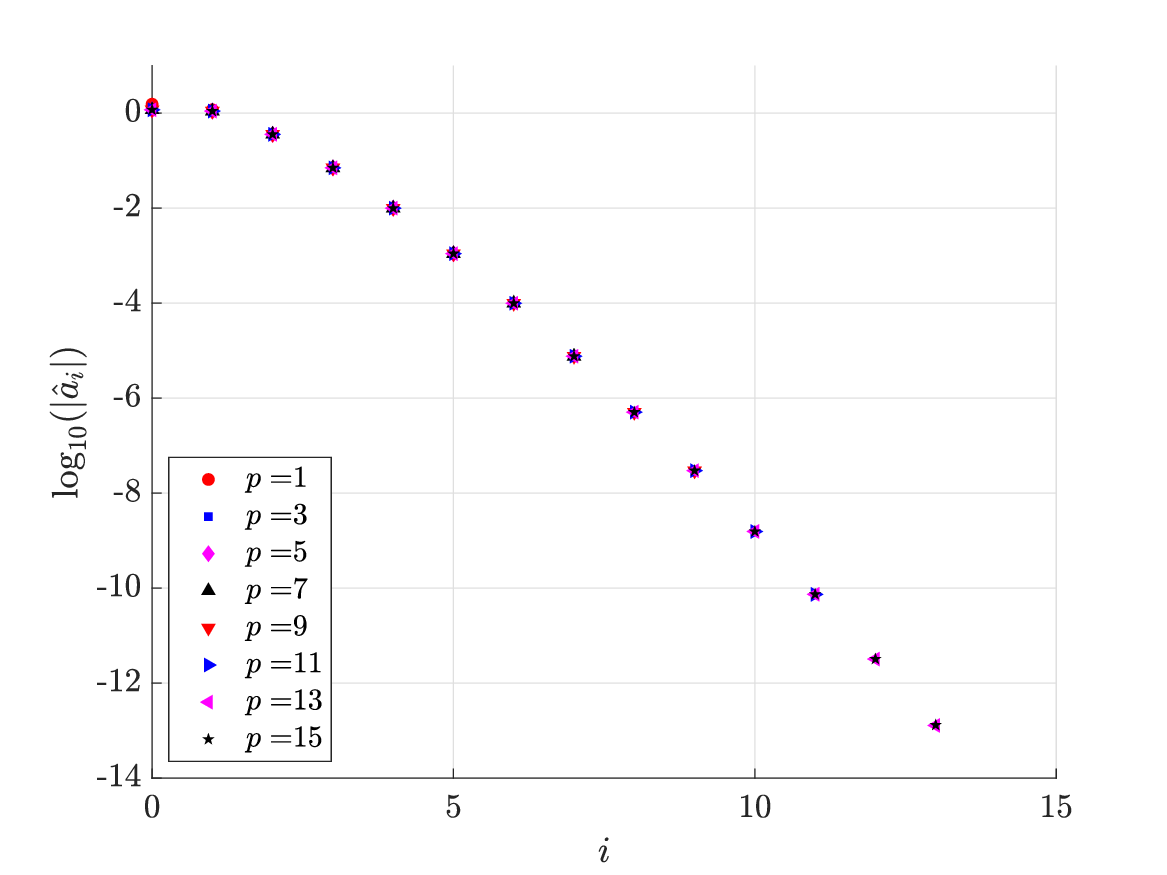}}
    \hfill
    \subfloat[$|\xi|$ \label{fig:coeffs2}]{\includegraphics[scale=0.4]{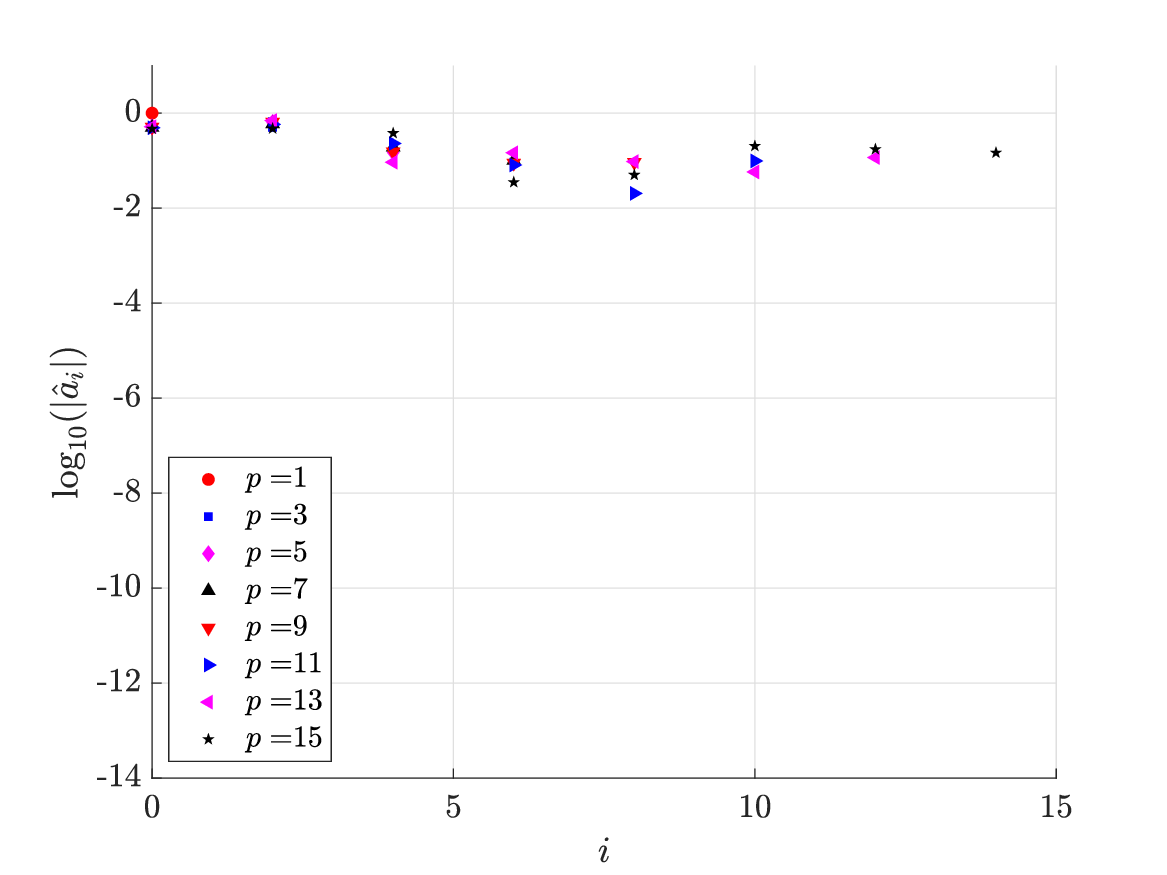}}
    \caption{Legendre polynomial coefficients of approximations of functions given in Eq.~\eqref{eq:funcs}.\label{fig:legcoeff}}
\end{figure*}

\begin{figure}
    \centering
    \includegraphics[scale = 0.4]{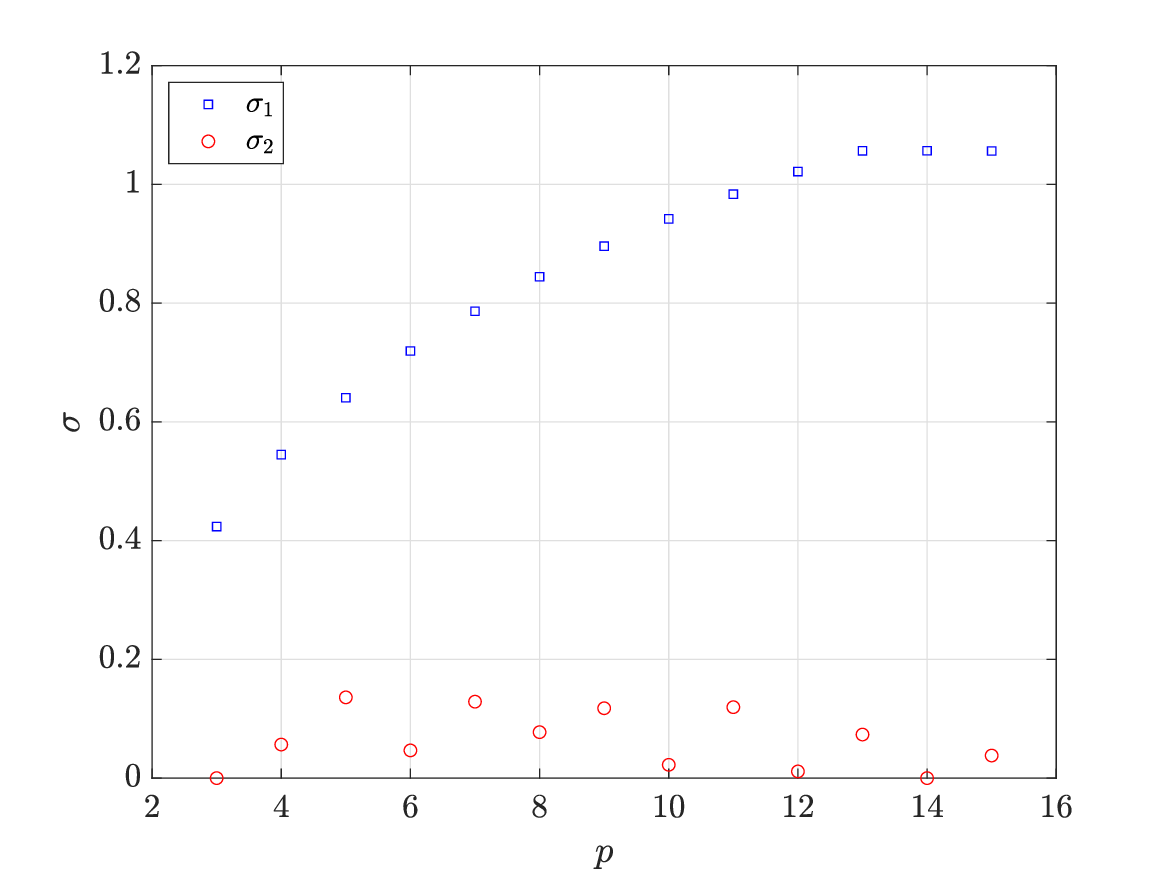}
    \caption{Decay rates as a function of the element order for the functions provided in Eq.~\eqref{eq:funcs}.}
    \label{fig:decayrate}
\end{figure}

To demonstrate the estimation procedure in practice, we consider the two following functions on $\xi \in [-1,+1]$:
\begin{equation}
    y_1(x) = \mathrm{exp}(\xi),\quad y_2(\xi) = |\xi|, \label{eq:funcs}
\end{equation}
where $y_1(\xi)$ and $y_2(\xi)$ are smooth and nonsmooth (particularly, not analytic) functions, respectively. A plot of the $\mathrm{log}_{10}(|\hat{a}_i|)$ for varying element degrees for both functions is provided in Fig.~\ref{fig:legcoeff}. 
Using an orthogonal least-squares fit to the coefficients in Fig.~\ref{fig:legcoeff}, we may then compute an estimate for the decay rate of the coefficients. The estimated decay rates as a function of the degree of each element for the functions in Eq.~\eqref{eq:funcs} are shown in Fig.~\ref{fig:decayrate}.
Note, in Fig.~\ref{fig:decayrate}, the decay rates are larger as a function of the element degree for the smoother function, $y_1(\xi)$. As it has been assumed that the solution is analytic, for some functions that are \textit{not} analytic such as $y_2(\xi)$, some of the estimated decay rates can be either negative or zero. In each of these instances, the decay rate is specified to be zero (which enforces $h-$refinement) and the function is properly described as nonsmooth. 

\subsection{Method for Mesh Refinement}
Suppose a solution has been produced on mesh $\bar{M}$ to the OCP in Eqs.~\eqref{eq:objective}-\eqref{eq:initialcondition} and the error has been computed on each element by Eq.~\eqref{eq:relerrorspace}. If the maximum relative error in any element is greater than a user-specified tolerance, $\epsilon_x$, then the element is modified by either increasing the polynomial degree or dividing the element into multiple smaller elements. In this case, we will use the decay rates of the Legendre polynomial coefficients to decide which is performed. The refinement technique is extended from Ref.~\cite{LiuRao2017} for use in the finite element discretization.

\subsubsection{Method for $p-$Refinement}
If the relative error in element $k$ is greater than the mesh error tolerance and the decay rates across all time points are greater than $\bar{\sigma}$, the solution is regarded as smooth in element $k$ and, if possible, the polynomial approximation used on mesh $\bar{M}+1$ is increased to reduce the solution error. Using the upper bound of the error in Eq.~\eqref{eq:errorlegendre}, on mesh $\bar{M}$, we have 
\begin{equation}
    \hat{e} = \dfrac{c10^{-\sigma(p^{(k,\bar{M})}+1)}}{\sqrt{1-10^{-2\sigma}}}. \label{eq:errorlegendre2}
\end{equation}
On the ensuing mesh, $\bar{M}+1$, it is desired to achieve a relative error, $\epsilon$:
\begin{equation}
    \epsilon = \dfrac{c10^{-\sigma(p^{(k,\bar{M}+1)}+1)}}{\sqrt{1-10^{-2\sigma}}}. \label{eq:epseqn}
\end{equation}
Equations \eqref{eq:errorlegendre2} and \eqref{eq:epseqn} can be solved for $p^{(k,\bar{M}+1)}$ to yield 
\begin{equation}
    p^{(k,\bar{M}+1)} = p^{(k,\bar{M})} +\dfrac{\mathrm{log}_{10}\left(\frac{\hat{e}}{\epsilon} \right)}{\sigma}.
\end{equation}
Replacing $\hat{e}$ with $\eta^{(k)}_{x}$ from our estimate in Eq.~\eqref{eq:relerrorspace} and taking the minimum decay rate over the temporal points provides 
\begin{equation}
    p^{(k,\bar{M}+1)} = p^{(k,\bar{M})} +\dfrac{\mathrm{log}_{10}\left(\frac{\eta^{(k)}_{x}}{\epsilon} \right)}{\min_{i=0}^{N_t}[\sigma_i]}, \label{eq:colpredicted}
\end{equation}
and to ensure that  $p^{(k,\bar{M}+1)} > p^{(k,\bar{M})}$, we add 
\begin{equation}
    p^{(k,\bar{M}+1)} = p^{(k,\bar{M})} + \left\lceil\dfrac{\mathrm{log}_{10}\left(\frac{\eta^{(k)}_{x}}{\epsilon} \right)}{\min_{i=0}^{N_t}[\sigma_i]} \right\rceil. \label{eq:colnew}
\end{equation}
If $p^{(k,\bar{M}+1)}$ is greater than a user-specified maximum, $p_{\max}$, the element is split into multiple elements with each containing the original number of nodes $p^{(k,\bar{M})} +1$ such that the sum of the polynomial degrees over the newly created subelements is greater than or equal to $p^{(k,\bar{M}+1)}$.

\subsubsection{Method for $h-$Refinement}
If, on the other hand, the relative error in element $T$ is greater than the mesh error tolerance and the decay rates across all time points are \textit{not} greater than $\bar{\sigma}$, the solution is regarded as nonsmooth in element $T$, and the element is broken into multiple smaller elements. To determine the number of elements to be added, the sum of the number of collocation points in the newly created elements should be greater than or equal to the polynomial degree predicted in Eq.~\eqref{eq:colpredicted}. However, to avoid constructing an infeasibly large mesh for nonsmooth problems, the threshold $\bar{\sigma}$ is used in place of $\min_{i =0}^{N_t}[\sigma_i]$. The sum is obtained by 
\begin{equation}
    \tilde{p} = p^{(k,\bar{M})} +\dfrac{\mathrm{log}_{10}\left(\frac{\eta^{(k)}_{x}}{\epsilon} \right)}{\bar{\sigma}}. \label{eq:addedmeshints}
\end{equation}
Each subelement that is created will contain the same number of collocation points as element $k$ on mesh $\bar{M}$. Thus, the number of newly created subelements, $H^{(k)}$, into which element $k$ is divided is computed as 
\begin{equation}
    H^{(k)} = \left\lceil \frac{\tilde{p}}{p^{(k,\bar{M})}} \right\rceil.
\end{equation}

\subsection{Method for Mesh Reduction}
In addition to a method to increase the degree of the state approximation and the number of elements, the work of Ref.~\cite{LiuHager2015} is used to decrease the size of the mesh in overcollocated regions through polynomial degree reduction and element reduction. For more information on the reduction process, the reader may refer to Appendix \ref{app:Reduction}.

\subsection{The Refinement/Reduction Algorithm}
\begin{figure*}[t]
    \centering
    \includegraphics[scale = 0.75]{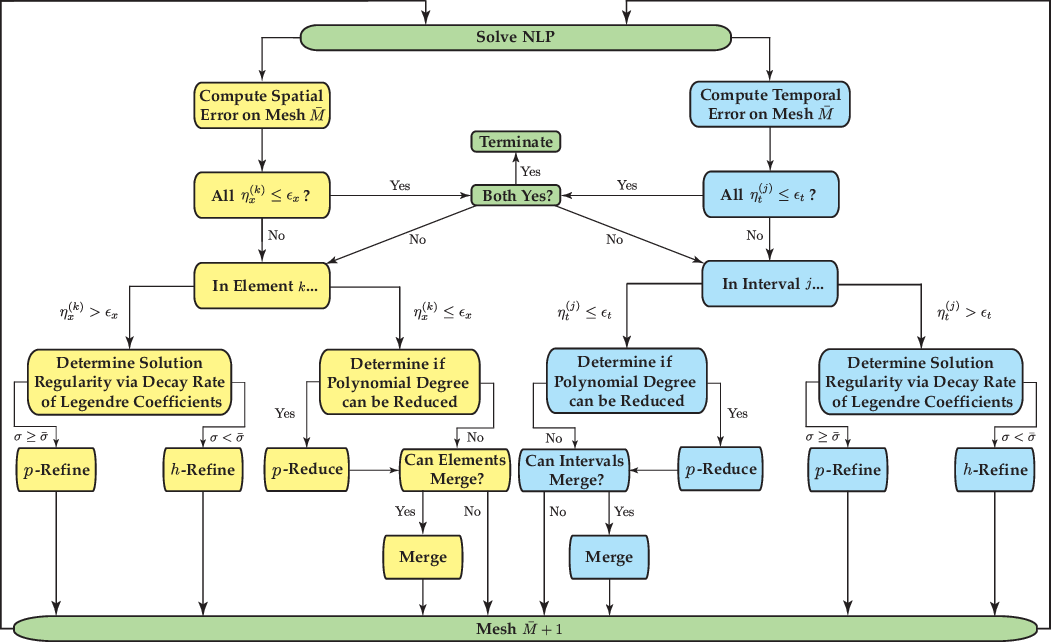}
    \caption{The refinement and reduction algorithm.}
    \label{fig:flowchart}
\end{figure*}

The refinement/reduction algorithm is described as follows. The NLP from Section \ref{sec:NLPsec} is solved using a well-developed NLP software such as \textit{IPOPT}, \textit{SNOPT}, or \textit{KNITRO}. Following the solution of the NLP, the error in both the temporal and spatial dimensions is computed via the methods described in Section \ref{sec:error}. Due to the reliance upon parallelization to evaluate each error estimate, it is generally more computationally efficient to evaluate the errors in each dimension sequentially in serial. That is, first, the error in the spatial dimension is computed, and then the error in the temporal error is computed, or vice versa, rather than both estimates being computed simultaneously on multiple cores. It is generally good practice to avoid nested parallelization, as evaluating each \textit{dimension} in serial allows all computational threads to be utilized for each dimension rather than dividing threads for simultaneous evaluation of the error in both dimensions. Unless the user has the computational resources to divide the cores effectively, nested parallelization will prove more inefficient than evaluating the error in each dimension sequentially. Once the estimate has been computed, if all error indicators in each element and interval are below the requested tolerance, the refinement algorithm will terminate. If not, the algorithm will continue with the refinement and reduction methods described in Section \ref{sec:refine}. In the instance that in one dimension the requested error tolerance is satisfied but in the other it is not, the dimension that satisfies the error tolerance will still be checked for potential mesh reduction. A flow chart outlining the refinement and reduction algorithm is provided in Fig.~\ref{fig:flowchart}. To avoid runaway refinement in the instance of an error asymptote, the user may supply a maximum number of refinement iterations to alternatively terminate the algorithm.     

\section{Numerical Examples}\label{sec:examples}

In this section, two numerical examples are shown that demonstrate the presented method in practice. The first problem analyzes the optimal control of the viscous Burgers' equation, while the second problem focuses on the optimal control of a nonlinear heat equation.
All numerical solutions were computed using MATLAB R2023b on a 2023 Apple M2 Ultra Mac Studio equipped with 24 cores and the open source NLP solver {\em IPOPT} with the linear solver MA57 and a limited-memory Hessian approximation. For all numerical examples presented, the regularity constant $\bar{\sigma} = 0.5$ in both the spatial and temporal dimensions, and the requested error tolerance $\epsilon_t = \epsilon_x = \epsilon$ is identical in both dimensions. 

\subsection{Viscous Burgers' Equation}

The first example analyzes the optimal boundary control of the viscous Burgers' equation. The problem was presented in Ref.~\cite{BuskensGriesse2006} and additionally solved in Ref.~\cite{Betts2020}. The goal is to minimize 
\begin{equation}
\mathcal{J} = \frac{1}{2}\int_0^1 \int_0^1 [y(x,t)-0.035]^2 \:\mathrm{d}x\mathrm{d}t + \frac{\gamma}{2} \int_0^1 [u_1^2(t)+u_2^2(t)]\:\mathrm{d}t,
\end{equation}
subject to the PDE
\begin{equation}
\partial_t y(x,t) = \nu \partial_x(\partial_x y(x,t)) - y(x,t)\partial_x y(x,t),
\end{equation}
with the initial conditions
\begin{equation}
y(x,0) = x^2(1-x)^2,
\end{equation}
the Neumann boundary conditions
\begin{align}
\partial_x y(0,t) &= u_1(t), \\
\partial_x y(1,t) & = u_2(t),
\end{align}
and control bounds for $i = 1,\:2$:
\begin{equation}
u_{\text{min}} \leq u_i(t) \leq u_{\text{max}}.
\end{equation}
The problem is defined on the domain 
\begin{equation}
\Omega = \{(x,t)\:\:|\:\:0\leq x \leq 1;\:\: 0\leq t \leq 1\},
\end{equation}
and the following parameters complete the definition:
$$
\gamma = 0.01, \:\:\:\:\:\nu = 0.1,\:\:\:\:\:u_{\text{max}} = 0.015,\:\:\:\:\:u_{\text{min}} = -0.015.
$$
The problem is solved on an initial mesh with 2 intervals and 6 collocation points per interval and 9 quadratic elements. The NLP tolerance was set to $10^{-12}$ with an acceptable tolerance of $10^{-10}$. The problem was solved to five requested relative error tolerances of $\epsilon = 10^{-4}, 10^{-5}, 10^{-6}, 10^{-7},\:\mathrm{and}\:10^{-8}$. 

For this problem, we compare the ``Local $hp$'' method presented in this paper with three other global refinement algorithms to demonstrate the effectiveness of the method in practice. The three global refinement algorithms are as follows. The ``Global $h$'' method is an algorithm which doubles the number of mesh intervals/elements in a given dimension if the maximum error over the intervals/elements is greater than the requested tolerance. If the maximum error in a particular dimension is less than the requested error tolerance, no intervals/elements in that dimension are added. The ``Global $p$'' method is a refinement algorithm which fixes the number of mesh intervals/elements in both dimensions but increases the degree of the state approximation in all intervals/elements by four if the maximum error is greater than the requested tolerance. If the maximum error in a particular dimension is less than the requested error tolerance, the polynomial degree in each interval/element is not increased in that dimension. Lastly, the ``Global $p$-$h$'' method combines the aspects of $p-$ and $h-$refinement to compare a global technique that harnesses both methods for convergence. If the maximum error in a given dimension is greater than the requested tolerance, the ``Global $p$-$h$'' method increases the degree of the state approximation in all intervals/elements by one until a specified maximum degree is reached. Once the specified maximum (for the degree of the state approximation) has been reached, and in the event the error tolerance is still not satisfied, the number of intervals/elements is doubled until the maximum error is less than the specified tolerance. If the maximum error in a particular dimension is less than the requested tolerance, neither the degree of the state approximation nor the number of intervals/elements are increased. In this work, the specified maximum for the degree of the state approximation in both dimensions is eight.    

The solution details are presented in Table~\ref{tab:burgerssoln}. The maximum error over all elements is denoted by $\eta_x^{\max}$, and the maximum error over all intervals is denoted by $\eta_t^{\max}$. The value of the optimal objective is denoted by $\mathcal{J}^*$. The number of mesh refinement iterations is denoted by $\bar{M}$. The CPU times presented in Table~\ref{tab:burgerssolnCPU} are averaged over 10 trials. The \textbf{adigator} column represents the total amount of time spent generating NLP derivative files during each trial. The \textbf{NLP} column indicates the total amount of time per trial spent solving the NLP. The \textbf{TR} and \textbf{SR} column indicates the total time spent refining the temporal and spatial mesh, respectively, including the evaluation of the error estimates, and lastly, the \textbf{Total} column indicates the average time to complete each trial. 

\begin{table*}[htbp]

\caption{\label{tab:burgerssoln}Solution details for the Burgers' equation tracking problem} 
\centering
\begin{tabular}{c|c|cccccccc}
\hline
Refinement & $\epsilon$ & $\bar{M}$ & $\eta_t^{\max}$ & $\eta_x^{\max}$ & $\mathcal{J}^*$ & $N_t$ & $J$ & $N_x$ & $K$ \\
\hline
None & - &  0 & $5.36 \times 10^{-5}$ & $4.43 \times 10^{-4}$ & $2.8940597 \times 10^{-5}$ & 12 & $2$ & 19 & 9 \\
\hline
Global $h$ &\multirow{4}{*}{$10^{-4}$} & 2 & $5.37 \times 10^{-5}$ & $2.08 \times 10^{-5}$ & $2.8969888 \times 10^{-5}$ & 12 & 2 & 73 & 36  \\

Global $p$ &  & 1 & $5.37 \times 10^{-5}$ & $4.95 \times 10^{-9}$ & $2.8970004 \times 10^{-5}$ & 12 & 2 & 55 & 9  \\

Global $p$-$h$ &  & 1 & $5.38 \times 10^{-5}$ & $1.52 \times 10^{-5}$ & $2.8969606 \times 10^{-5}$ & 12 & 2 & 28 & 9  \\

Local $hp$ &  & 2 & $5.37 \times 10^{-5}$ & $8.93 \times 10^{-5}$ & $2.8970043 \times 10^{-5}$ & 12 & 2 & 31 & 12  \\
\hline
Global $h$ &\multirow{4}{*}{$10^{-5}$} & 3 & $8.70 \times 10^{-6}$ & $4.72 \times 10^{-6}$ & $2.8969376 \times 10^{-5}$ & 48 & 8 & 145 & 72  \\

Global $p$ &  & 3 & $6.54 \times 10^{-6}$ & $3.83 \times 10^{-8}$ & $2.8969341 \times 10^{-5}$ & 36 & 2 & 55 & 9  \\

Global $p$-$h$ &  & 4 & $5.35 \times 10^{-6}$ & $2.54 \times 10^{-6}$ & $2.8969375 \times 10^{-5}$ & 64 & 8 & 37 & 9  \\

Local $hp$ & & 2 & $8.69 \times 10^{-6}$ & $8.45 \times 10^{-6}$ & $2.8969342 \times 10^{-5}$ & 24 & 4 & 51 & 15  \\
\hline
Global $h$ & \multirow{4}{*}{$10^{-6}$} & \multicolumn{8}{c}{----------------- No solution found, memory exceeded -----------------}  \\

Global $p$ &  & 12 & $9.31 \times 10^{-7}$ & $4.96 \times 10^{-7}$ & $2.8969388 \times 10^{-5}$ & 108 & 2 & 55 & 9 \\

Global $p$-$h$ &   & 7 & $3.91 \times 10^{-7}$ & $8.67 \times 10^{-7}$ & $2.8969490 \times 10^{-5}$ & 512 & 64 & 46 & 9  \\

Local $hp$ &   & 5 & $8.72 \times 10^{-7}$ & $7.87 \times 10^{-7}$ & $2.8969338 \times 10^{-5}$ & 48 & 8 & 77 & 21  \\
\hline
Global $h$ & \multirow{4}{*}{$10^{-7}$}  & \multicolumn{8}{c}{----------------- No solution found, memory exceeded -----------------}  \\

Global $p$ &  & \multicolumn{8}{c}{-------------------- No solution found, $IPOPT$ failure -----------------} \\

Global $p$-$h$ &  & \multicolumn{8}{c}{----------------- No solution found, memory exceeded -----------------}  \\

Local $hp$ &   & 8 & $8.14 \times 10^{-8}$ & $8.17 \times 10^{-8}$ & $2.8969407 \times 10^{-5}$ & 82 & 13 & 111 & 27  \\
\hline
Global $h$ & \multirow{4}{*}{$10^{-8}$}  & \multicolumn{8}{c}{----------------- No solution found, memory exceeded -----------------}  \\

Global $p$ &   & \multicolumn{8}{c}{-------------------- No solution found, $IPOPT$ failure -----------------} \\

Global $p$-$h$ &   & \multicolumn{8}{c}{----------------- No solution found, memory exceeded -----------------}  \\

Local $hp$ &  &  10 & $5.23 \times 10^{-9}$ & $9.37 \times 10^{-9}$ & $2.8969446 \times 10^{-5}$ & 131 & 19 & 168 & 38  \\
\hline
\end{tabular}
\end{table*}

\begin{table*}[h!]

\caption{\label{tab:burgerssolnCPU}CPU details for global vs.~local refinement techniques} 
\centering
\begin{tabular}{c|c|ccccc}
\hline

Refinement & $\epsilon$ & \textbf{adigator} (s) & \textbf{NLP} (s) & \textbf{TR} (s) & \textbf{SR} (s) & \textbf{Total} (s) \\
\hline 
None & - & 2.9045 & 0.1103 & - & - & 3.0153 \\
\hline
Global $h$ & \multirow{4}{*}{$10^{-4}$} & 8.3083 & 0.4129 & 0.5961 & 0.5983 & 10.013 \\

Global $p$ & & 5.6157 & 0.3226 & 0.3056 & 0.3024 & 6.7020 \\

Global $p$-$h$ &  & 5.7119 & 0.2139 & 0.1965 & 0.2666 & 6.4306 \\

Local $hp$ &  & 8.3465 & 0.3355 & 0.2965 & 0.3968 &  9.4341 \\
\hline
Global $h$ & \multirow{4}{*}{$10^{-5}$} & 12.143 & 8.7787 & 3.1390 & 2.7657 & 27.094 \\

Global $p$ &  & 11.068 & 3.3714 & 3.0918 & 1.0859 & 18.980 \\

Global $p$-$h$ &  & 13.819 & 1.5385 & 0.7252 & 1.3712 & 17.606 \\

Local $hp$ &  & 8.3693 & 2.0365 & 0.4468 & 0.5626 & 11.529 \\
\hline
Global $h$ & \multirow{4}{*}{$10^{-6}$} & \multicolumn{5}{c}{------------------------ No solution found ------------------------}\\

Global $p$ &  & 51.665 & 410.88 & 129.22 & 14.335 & 607.35 \\

Global $p$-$h$ &  & 39.413 & 89.650 & 2.8821 & 123.67 & 256.06 \\

Local $hp$ &  & 17.671 & 7.0881 & 2.5911 & 2.8557 & 30.701 \\
\hline
Global $h$ & \multirow{4}{*}{$10^{-7}$} & \multicolumn{5}{c}{------------------------ No solution found ------------------------}\\

Global $p$ &  & \multicolumn{5}{c}{------------------------ No solution found ------------------------} \\

Global $p$-$h$ &  & \multicolumn{5}{c}{------------------------ No solution found ------------------------} \\

Local $hp$ &  & 32.926 & 40.917 & 14.340 & 12.301 & 101.992 \\
\hline
Global $h$ & \multirow{4}{*}{$10^{-8}$} & \multicolumn{5}{c}{------------------------ No solution found ------------------------}\\

Global $p$ &  & \multicolumn{5}{c}{------------------------ No solution found ------------------------} \\

Global $p$-$h$ &  & \multicolumn{5}{c}{------------------------ No solution found ------------------------} \\

Local $hp$ &  & 89.3889 & 1005.3 & 112.04 & 50.090 & 1260.5 \\
\hline
\end{tabular}
\end{table*}

\begin{figure*}[h!]
    \subfloat[Optimal state \label{fig:burgersstate}]{\includegraphics[scale=0.4]{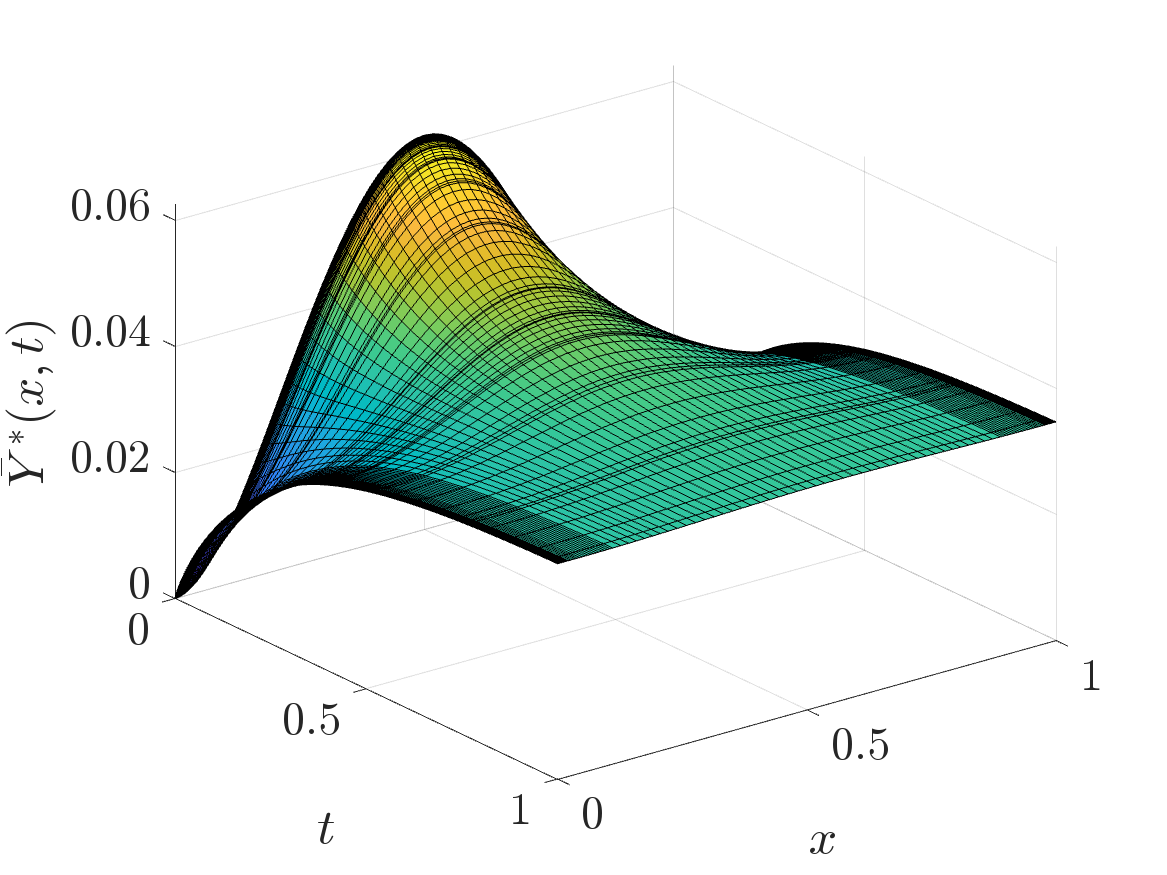}}
    \hfill
    \subfloat[Optimal controls\label{fig:burgerscontrol}]{\includegraphics[scale=0.4]{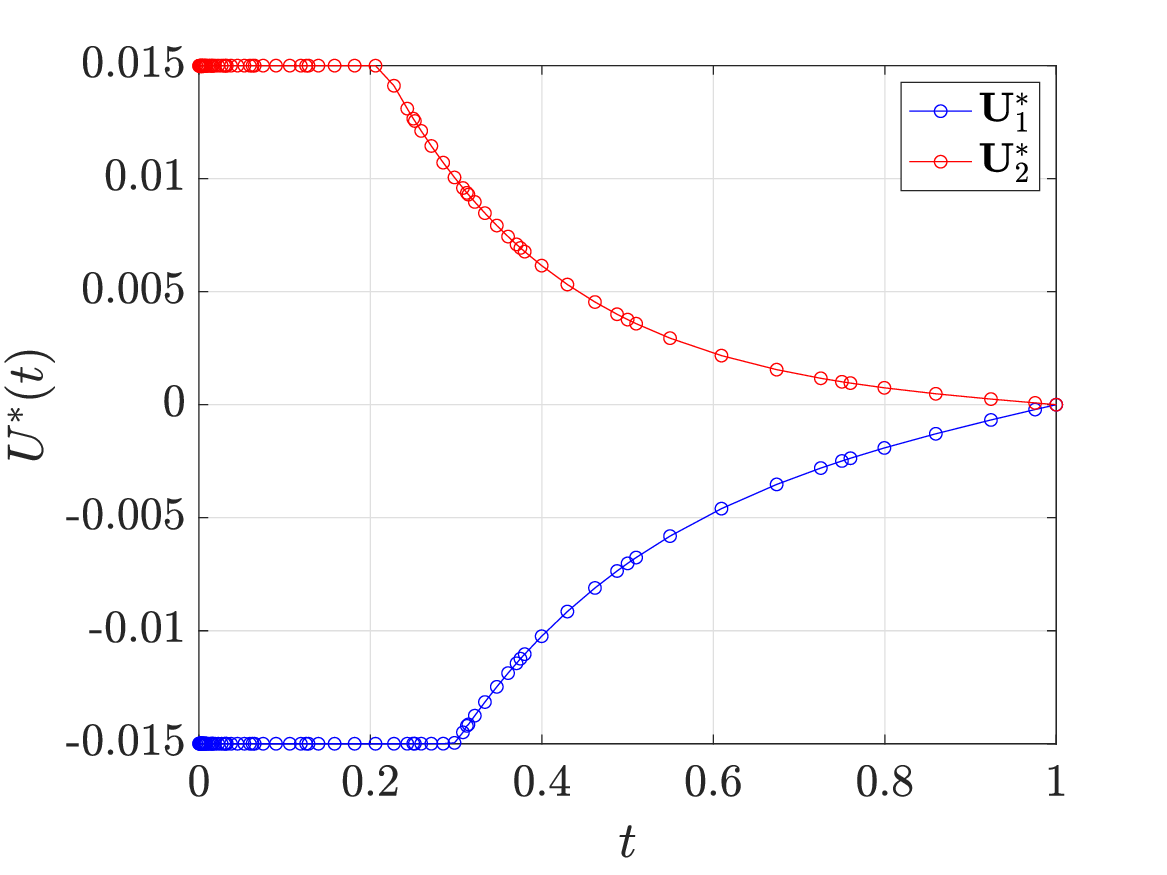}}
    \caption{The optimal state and optimal control for the Burgers' equation tracking problem with $\epsilon = 10^{-7}$. \label{fig:BurgersSoln}}
\end{figure*}

\begin{figure*}[h!]
    \centering
    \includegraphics[scale = 0.675]{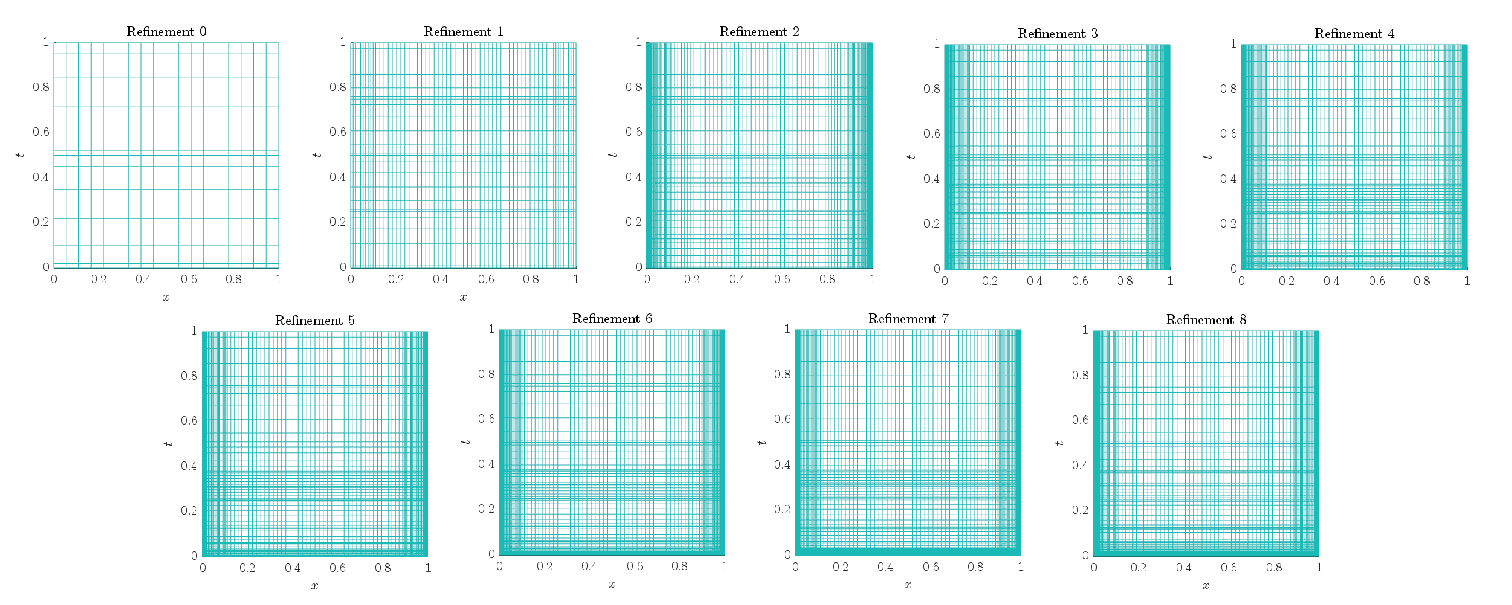}
    \caption{Refinement sequence for $\epsilon = 10^{-7}$ for the Burgers' equation tracking problem.}
    \label{fig:burgersrefine}
\end{figure*}

A plot of the ``Local $hp$'' method result for $\epsilon = 10^{-7}$ is provided in Fig.~\ref{fig:BurgersSoln}. A plot of the refined mesh is provided in Fig.~\ref{fig:burgersrefine}. As is evident from Table \ref{tab:burgerssoln}, as the requested error tolerance decreases, the need for local refinement techniques increases. For the cases in which $\epsilon = 10^{-7}, \:10^{-8}$, all global algorithms fail to produce a solution to the problem due to the problem size. Interestingly, the ``Global $h$'' and ``Global $p$-$h$'' methods fail as a result of the NLP derivative file size (\textit{typically} due to the number of NLP variables present in the problem), which exceeds MATLAB's ``safe operation'' limit as specified by \textbf{adigator}; whereas, the ``Global $p$'' method fails as a result of the total number of non-zero elements in the constraint Jacobian of the NLP. That is, the ``Global $p$'' method severely decreases the sparsity of the NLP as the mesh grows, which hinders the ability of $IPOPT$ to find a viable solution. In the event that the density and size of the NLP becomes \textit{too} large, $IPOPT$ fails altogether, which is the case for $\epsilon \leq 10^{-7}$.  

The interdimensional error dependency is made clear by the results of the ``Global $p$'' method in Table \ref{tab:burgerssoln} as the error tolerance decreases from $10^{-4}$ to $10^{-6}$. In all refinement results from $10^{-4}$ to $10^{-6}$, despite having the same spatial mesh size and characteristics on each converged mesh, the maximum value of the spatial error increases from $4.9465 \times 10^{-9}$ for $\epsilon = 10^{-4}$ to $4.9638 \times 10^{-7}$ for $\epsilon = 10^{-6}$. The addition of more temporal points on the converged mesh provides additional locations to evaluate the error in the spatial dimension. In general, the addition of these points can increase the error in the spatial dimension if the temporal points are added in regions of the domain where the solution is nonregular. Thus, it becomes clear that there can exist multiple avenues to achieve mesh sizes that satisfy the requested error tolerance. The ``Global $p$'' and ``Global $p$-$h$'' methods do not leverage the underlying solution structure to inform the refinement process. By blindly increasing the mesh size in dimensions that do not satisfy the requested tolerance, the converged mesh can take on the structure of those in Table \ref{tab:burgerssoln} for $\epsilon = 10^{-6}$. In order to provide enough points on the grid to resolve the solution in structured regions, refinement in one dimension may be exhausted to fulfill the total error requirement (i.e.~both dimensions satisfying the requested tolerance). When refining the mesh based on the underlying solution structure, such as the ``Local $hp$'' method presented in this work, points may be more astutely placed in both dimensions to mitigate the total size and difficulty of the problem. 

Another benefit of the method presented in this work is the ability to leverage sparsity in the NLP to more efficiently produce solutions to the problem. A downfall of the ``Global $p$'' approach is that the density of the NLP constraint Jacobian grows as more points are added. Thus, even if the total number of points on the converged mesh for the ``Global $p$'' method is lower than other converged mesh sizes that use $h-$ and $p-$refinement, the solution can be more computationally expensive to produce. This phenomenon is indicated by the results for the ``Global $p$'' and ``Global $p$-$h$'' methods for $\epsilon = 10^{-6}$ in Tables \ref{tab:burgerssoln} and \ref{tab:burgerssolnCPU}, where as there are $23,598$ total points on the converged mesh using the ``Global $p$-$h$'' method, there are only $5,995$ points on the converged mesh using the ``Global $p$'' method. However, on average over ten trials of the method, the ``Global $p$'' method took nearly 350 seconds longer to find a solution. The contrast between the two total CPU times is most accredited to the disparity in the time required to the solve the NLP. As the density of the NLP for the ``Global $p$'' method is much larger than the density of the NLP for the ``Global $p$-$h$'' method, the NLP becomes more computationally expensive to solve. The benefit of the method produced in this paper is that the underlying solution structure may be leveraged to mitigate the overall problem size and resulting sparsity of the NLP to produce accurate solutions to PDE-constrained OCPs in a computationally efficient procedure.     


In Table~\ref{tab:burgerssolnCPU}, it is clear that there is a significant computational cost to refine the mesh beyond $\epsilon = 10^{-7}$. The cost grows as a result of the distribution of the error on the mesh. Adding more points in regions that are already significantly refined can increase the difficulty of the NLP solution, which significantly increases the cost of solving the problem. Though the cost of evaluating the estimate grows as a function of the mesh size, the parallel evaluation of the error estimates aids to mitigate the overall computational cost of the ``Local $hp$'' method, which significantly outperforms global techniques (particularly for high-accuracy solutions) even if the cost of evaluating the estimates for the global procedures are ignored.  



As is evident from both Figs.~\ref{fig:BurgersSoln} and \ref{fig:burgersrefine}, the regions in which the mesh is most refined is near the initial time and the near and far spatial boundaries. The solution exhibits the most structure in these regions, which yields expected refinement. In addition to capturing the change in state structure, the refinement algorithm effectively captures the regions in which the control leaves the constraint bound around Refinement 4 in Fig.~\ref{fig:burgersrefine}. Though the error is technically defined as a ``state'' error, the means in which it is computed relies upon computation of the residual in the differential equations at points that are distinct from the NLP. The residual equations depend upon an interpolation of the control variable; thus, poor interpolation of the control (i.e.~poor \textit{resolution} of the control) can lead to greater values of the residual error, which induces refinement. In more static regions, the error is quickly resolved to the requested tolerance. As is evident from the solution in Fig.~\ref{fig:BurgersSoln}, the state approaches steady state at $t > 0.5$. From Fig.~\ref{fig:burgersrefine}, it is apparent that the mesh tolerance in this region is satisfied after the first refinement iteration, and little to no future refinement is required to resolve the error. For a comparison to results with refinement in one dimension, we refer the reader to Refs.~\cite{Betts2020, DaviesPollock2026}.

\subsection{Nonlinear Heat Equation}\label{sec:heat1}
The second example analyzes the optimal boundary control of a nonlinear heat equation. The problem is solved in Refs.~\cite{Heinkenschloss1996,Betts2020}. It may be viewed as a simplified model for the heating of a probe in a kiln. The goal is to minimize the deviation from a desired temperature profile as defined by the objective 
\begin{equation}
    \mathcal{J} = \frac{1}{2}\int_0^{t_f} \left\{[y(1,t)-y_d(t)]^2+\gamma u^2(t)\right\}\:\mathrm{d}t,
\end{equation}
by choosing the control function between allowable limits
\begin{equation}
    u_{\min} \leq u(t) \leq u_{\max},
\end{equation}
that satisfies the nonlinear heat equation
\begin{equation}
        q(x,t) = (a_1+a_2y)\partial_ty-a_3\partial_x(\partial_xy) 
        -a_4\left(\partial_xy\right)^2-a_4y\partial_x(\partial_xy), \label{eq:heatpde}
\end{equation}
and the initial and boundary conditions 
\begin{align}
        (a_3+a_4y)\left. \partial_xy\right|_{x=0} &= g[y(0,t)-u(t)], \\
    (a_3+a_4y)\left. \partial_xy\right|_{x=1} &= 0, \\
    y(x,0) &= y_{I}(x),
\end{align}
where the following definitions: 
    \begin{align*}
        y_d(t) =&\:2-\mathrm{e}^{\rho t}, \\
        y_I(x) =&\:2 + \mathrm{cos}(\pi x), \\
        q(x,t) =&\:[\rho(a_1+2a_2)+\pi^2(a_3+2a_4)]\mathrm{e}^{\rho t}\mathrm{cos}(\pi x) -a_4\pi^2\mathrm{e}^{2\rho t}+(2a_4\pi^2+\rho a_2)\mathrm{e}^{2\rho t}\mathrm{cos}^2(\pi x),
    \end{align*}
    and parameters in Eq.~\eqref{eq:neccparameters} complete the problem statement:
\begin{equation}
    \begin{array}{lclclcl}
        a_1 & = & 4 &, & \rho & = & -1, \\
        a_2 & = & 1 & , & t_f & = & 0.5, \\ 
        a_3 & = & 4 & , & \gamma & = & 10^{-3}, \\
        a_4 & = & -1 & , & g & = & 1, \\ 
        u_{\min} & = & -\infty & , & u_{\max} & = & 0.1.
    \end{array}
    \label{eq:neccparameters}
\end{equation}

\begin{table*}[h!]

\caption{\label{tab:heatsoln}Solution details for the nonlinear heat equation problem} 
\centering
\begin{tabular}{ccccccccc}
\hline
$\epsilon$ & $\bar{M}$ & $\eta_t^{\max}$ & $\eta_x^{\max}$ & $\mathcal{J}^*$ & $N_t$ & $J$ & $N_x$ & $K$ \\
\hline 
- & 0 & $4.1026 \times 10^{-4}$ & $0.0022$ & $3.8648480 \times 10^{-5}$ & 12 & $3$ & 19 & 9 \\

$10^{-5}$ & 5 & $6.9917 \times 10^{-6}$ & $7.0138 \times 10^{-6}$ & $3.8654831 \times 10^{-5}$ & 30 & 7 & 41 & 9  \\

$10^{-6}$ & 6 & $4.7932 \times 10^{-7}$ & $4.9043 \times 10^{-7}$ & $3.8654915 \times 10^{-5}$ & 49 & 10 & 57 & 10  \\

$10^{-7}$ & 8 & $9.7474 \times 10^{-8}$ & $1.2954 \times 10^{-8}$ & $3.8654934 \times 10^{-5}$ & 73 & 15 & 84 & 12  \\
\hline
\end{tabular}
\end{table*}

\begin{table*}[h!]

\caption{\label{tab:heatsolnCPU}CPU details for the nonlinear heat equation problem} 
\centering
\begin{tabular}{cccccc}
\hline
$\epsilon$ & \textbf{adigator} (s) & \textbf{NLP} (s) & \textbf{TR} (s) & \textbf{SR} (s) & \textbf{Total} (s) \\
\hline 
- & 2.3643 & 0.0705 & - & - & 2.4352 \\

$10^{-5}$ & 15.3558 & 5.8327 & 1.3949 & 2.2528 & 25.4405 \\

$10^{-6}$ & 18.6243 & 8.8572 & 4.7755 &  7.7987 & 41.8397 \\

$10^{-7}$ & 25.7815 & 41.1305 & 21.2660 & 35.1468 & 129.8608 \\
\hline
\end{tabular}
\end{table*}

\begin{figure*}[h!]
    \subfloat[Optimal state \label{fig:heatstate}]{\includegraphics[scale=0.4]{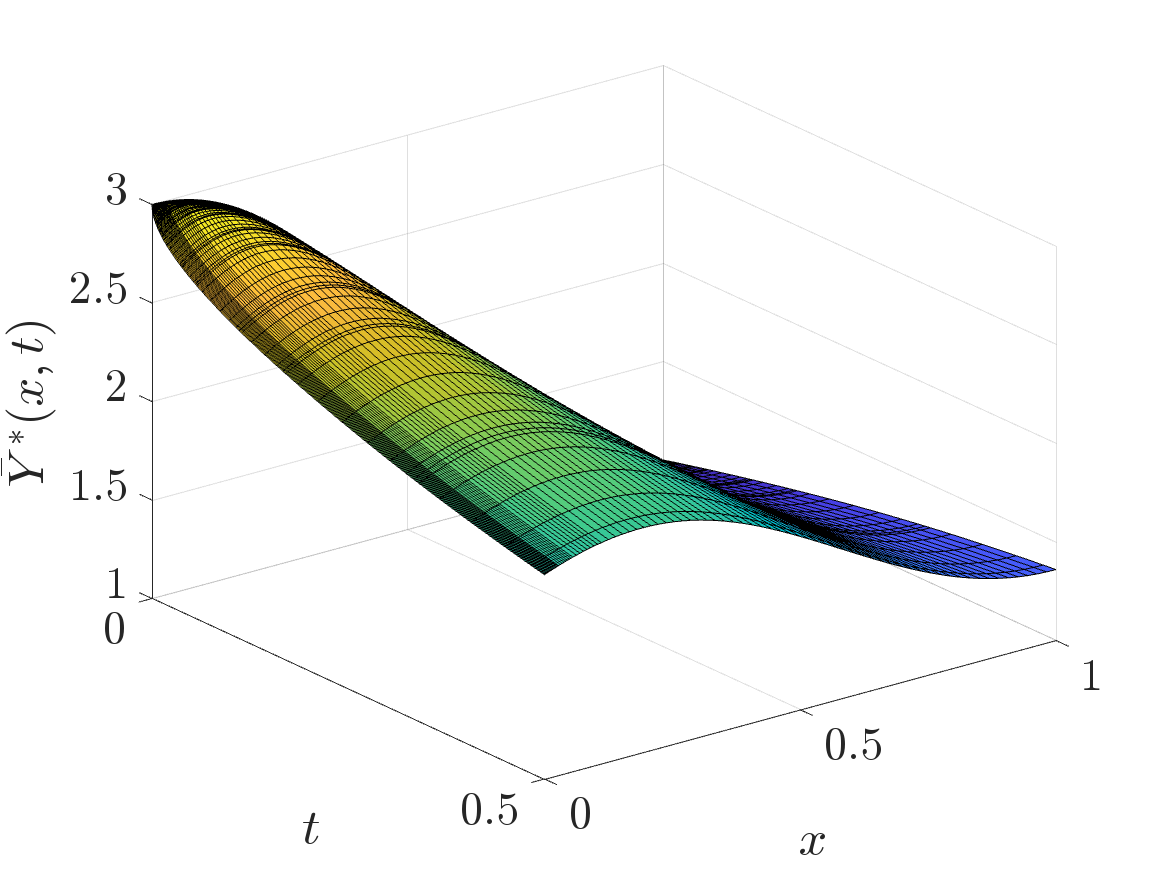}}
    \hfill
    \subfloat[Optimal control\label{fig:heatcontrol}]{\includegraphics[scale=0.4]{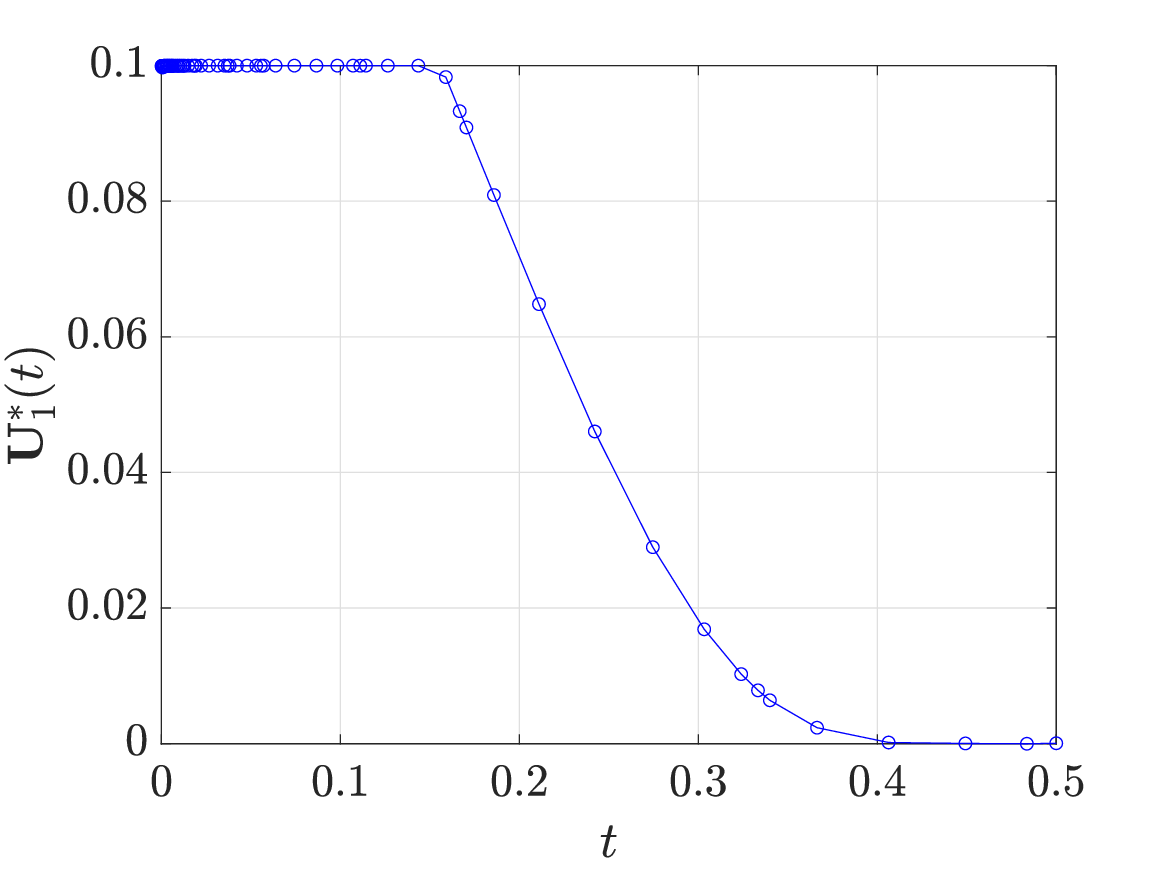}}
    \caption{The optimal state and optimal control for the heat equation with $\epsilon = 10^{-7}$. \label{fig:heateqn}}
\end{figure*}

\begin{figure*}[h!]
    \centering
    \includegraphics[scale = 0.65]{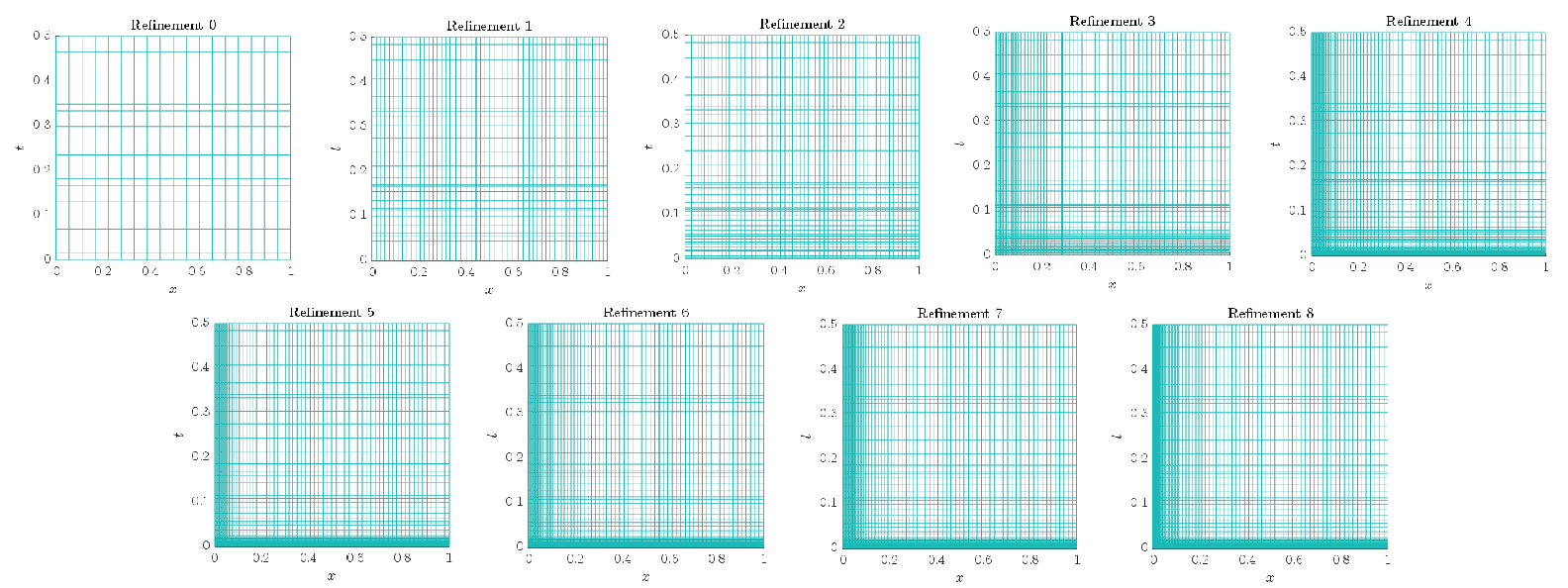}
    \caption{Refinement sequence for the nonlinear heat equation problem with $\epsilon = 10^{-7}$.}
    \label{fig:heatrefine}
\end{figure*}

The problem is solved with the ``Local $hp$'' method for requested mesh tolerances of $10^{-5},\:10^{-6},\:\mathrm{and}\:10^{-7}$ on an initial mesh with 3 mesh intervals containing 4 collocation points each and 9 quadratic elements. The requested NLP tolerance was set to $10^{-12}$ with an acceptable tolerance of $10^{-10}$. Solution details are provided in Table~\ref{tab:heatsoln} and \ref{tab:heatsolnCPU}. The refinement sequence for the solution generated in Fig.~\ref{fig:heateqn} is provided in Fig.~\ref{fig:heatrefine}. A common theme in the results of the adaptive algorithm presented here is that refinement is often performed in regions of the state that experience the most change. Often, this is near the boundaries where the control input is provided. As is evident from the sequence in Fig.~\ref{fig:heatrefine}, the majority of the refinement occurs near the initial time and the near spatial boundary. 
In less dynamic regions of the mesh, the refinement algorithm converges to a mesh size that satisfies the requested error tolerance within the first few mesh iterations. An example of this behavior occurs in the mesh region where $x \gtrsim 0.15$ and $t \gtrsim 0.05$, where the mesh refinement procedure is completed after just one to two iterations.

A portion of the novelty of the adaptive method presented here is that it can be shown to achieve accurate results with fewer temporal points in comparison to other methods, effectively reducing the problem size. Compared to the results in Ref.~\cite{Betts2020}, which does \textit{not} provide refinement in the spatial dimension, the adaptive method presented in this work is capable of producing a temporal mesh with fewer points when solving to similar tolerances while \textit{also} providing refinement in the spatial dimension.
A refinement algorithm in one dimension may lead to an improper estimation of the error due to poor resolution of the unrefined dimension. For an example of a refinement algorithm that only applies temporal adaptivity with the method presented in this work, see Ref.~\cite{DaviesPollock2026} for technical demonstration. Though the converged mesh is larger than the results presented in Ref.~\cite{DaviesPollock2026}, the error is reduced and equilibrated in both the temporal and spatial dimension, which provides assurances that the solution on the converged mesh is accurate in both dimensions. The results presented here demonstrate a mesh of reduced size while supplying refinement in \textit{both} dimensions, which has not been done previously to solve this problem. 

The CPU cost, in this case, is higher for the spatial refinement algorithm than it is for the temporal algorithm. This is due to the larger number of nonlinear terms in the PDE in Eq.~\eqref{eq:heatpde}. In computing the solution to the spatial residual problems, the nonlinear terms in the variational form are computed directly on each root-finding iteration. The benefit of such an approach is that a quantifiable assurance is provided that the solution to the linearized transformed system of the NLP converges to the solution of the fully nonlinear PDE (at \textit{least} to the extent of the requested error tolerance). The drawback is that the computational cost of determining the solution to each residual problem is greater. However, this approach reduces the cost of solving the NLP while also providing accuracy assurances. 

\section{Conclusion}\label{sec:conclusions}
In this work, the adaptive finite element method and adaptive orthogonal collocation approaches are unified in a single framework for PDE-constrained OCPs. The implementation of the scheme in a direct transcription approach avoids the derivation of the conditions for optimality, which are often complicated and specific for each OCP. When implemented in a direct transcription framework, the scheme provides a versatile approach for solving a wide range of parabolic PDE-constrained OCPs.

In the both dimensions, an implicit residual estimation method is used to compute an estimate for the error in each element. The estimate is novel within an orthogonal collocation framework. The implicit residual error estimation method can be used to compute a reliable estimate that does not suffer from the drawbacks of explicit estimators, which often involve extensive bounding and weighting techniques. Often, the implicit residual error estimation method has been avoided due to the computational expense of solving several local residual problems. In this work, the computational expense is mitigated through parallelization of the estimation procedure, which is shown to be more computationally efficient than global refinement procedures. Based on the computed error indicator, the mesh is subsequently refined based on an estimate of the regularity of the solution in the element. The regularity of the solution is estimated by the decay rate of the coefficients of a Legendre polynomial expansion of the solution. To the extent of the authors' knowledge, no work has combined the presented error estimate and refinement approaches in a direct transcription framework for OC for {\em any} set of differential equation constraints. The refinement algorithm presented in this work allows for an equilibration of the error in the temporal and spatial dimensions, which assures accuracy throughout the domain.

In two numerical examples, the $hp-$adaptive algorithm in this work is shown to reliably compute solutions that result in error reduction by up to nearly five orders of magnitude. A comparison to global refinement techniques is completed, and the computational performance of the presented method is discussed. Despite presenting refinement in both the temporal and spatial dimension, the effectiveness of the $hp-$adaptive orthogonal collocation method is demonstrated by computing accurate solutions with fewer points than existing, purely temporal refinement algorithms. Thus, the method presented in this paper offers a generalizable, efficient, and accurate framework for computing solutions to OCPs constrained by parabolic partial differential equations.




\section*{Funding Data}

\begin{itemize}
    \item U.S. Air Force Research Laboratory (Grant No.~FA8651-25-1-0002).
    \item U.S.~National Science Foundation (Grant No.~CMMI-2031213).
    \item U.S.~Office of Naval Research (Grant No.~N00014-22-1-2397).
\end{itemize}










 \appendix   

 \section{Method for Mesh Reduction \label{app:Reduction}}
Here, a method for reducing the mesh size in overcollocated regions is presented in the context of the spatial discretization. For a reference with application to the temporal dimension, refer to Ref.~\cite{LiuHager2015}.
\subsection{Reducing the Polynomial Degree}
In element $k$, if the accuracy tolerance has been satisfied, it is desirable to determine whether the tolerance may still be satisfied with a reduced degree of the state approximation. To determine whether the degree can be reduced, we rely upon a power series representation of the state in element $k$. The state approximation in element $k$ (at a particular point in time) is given by 
\begin{align}
    y_h^{(k)}(x) &= \sum_{i=1}^{p^{(k)}+1}\phi_i^{(k)}\left(\frac{x-x^{(k)}_1}{h^{(k)} } \right)Y_i^{(k)}, \label{eq:stateapproxorig} \\
    \phi_i^{(k)} &= \prod_{\substack{j=1 \\ i \neq j}}^{p^{(k)}+1} \frac{r - r_{j}^{(k)}}{r_{i}^{(k)}-r_{j}^{(k)}}.
\end{align}
Alternatively, the Lagrange polynomial $\phi^{(k)}_i(r)$ may be written as 
\begin{equation}
   \phi^{(k)}_i(r) = \sum_{l = 0}^{p^{(k)}}a_{li}r^l, \label{eq:stateapproxforreduce}
\end{equation}
where the coeffiecients $a_{li}$ depend solely on the element support points and can be computed by relation to a power series with identical roots to $\phi^{(k)}_i(r)$. That is, $Q_i(r)$ is a power series whose roots are $\left\{r^{(k)}_j\right\}_{\substack{j=1 \\ j \neq i}}^{p^{(k)}+1}$ and can be expressed as 
\begin{equation}
    Q_i(r) = \sum_{l = 0}^{p^{(k)}} Q_{li}r^l. \label{eq:powerseries}
\end{equation}
By the isolation property, we have that $\phi^{(k)}_i(r^{(k)}_i) = 1,$ so we can write
\begin{equation}
    \phi_i^{(k)}(r) = \frac{1}{Q_i(r_i)}Q_i(r) = \sum_{l =0}^{p^{(k)}}\frac{Q_{li}}{Q_i(r_i)}r^l,
\end{equation}
and by Eq.~\eqref{eq:stateapproxforreduce}, we have
\begin{equation}
    a_{li} = \frac{Q_{li}}{Q_i(r_i)}. \label{eq:coeffcompute}
\end{equation}
Combining Eqs.~\eqref{eq:stateapproxorig} and \eqref{eq:stateapproxforreduce} provides 
\begin{equation}
y_h^{(k)}(x) = \sum_{l = 0}^{p^{(k)}}b_{l}\left(\dfrac{x-x_1^{(k)}}{h^{(k)}} \right)^l, \quad b_l = \sum_{i=1}^{p^{(k)}+1}a_{li}Y_i^{(k)},
\label{eq:combinedstateapprox}
\end{equation}
where importantly, $(x-x_1^{(k)})/h^{(k)} \leq 1$ for $x \in \Omega_k$. If we partially expand Eq.~\eqref{eq:combinedstateapprox}, we obtain
\begin{equation}
   y_h^{(k)}(x) = \sum_{l = 0}^{p^{(k)}-1}b_{l}\left(\dfrac{x-x_1^{(k)}}{h^{(k)}} \right)^l + b_{p^{(k)}} \left(\dfrac{x-x_1^{(k)}}{h^{(k)}} \right)^{p^{(k)}}, \label{eq:partiallyexpanded}
\end{equation}
By removing the $p^{(k)}$ term from Eq.~\eqref{eq:partiallyexpanded}, the error in the state approximation of $y_h^{(k)}(x)$ is at most $b_{p^{(k)}}$. Note that this process can be repeated at all values of time to yield a vector of coefficients
\begin{equation}
\mathbf{b}_{p^{(k)}} = \begin{bmatrix}
    b_{p^{(k)}0} & \ldots & b_{p^{(k)}N_t}
\end{bmatrix}^T \in \mathbb{R}^{(N_t+1) \times 1}.
\end{equation}
To compare with the relative error tolerance, $\epsilon_x$, we normalize the coefficients $\mathbf{b}_{p^{(k)}}$ by $\mathbf{b}_{p^{(k)}} \odot \boldsymbol{\chi}_x^{(k)\odot-1}$ from Eq.~\eqref{eq:relparam}. Then, starting with the highest power, terms are removed from Eq.~\eqref{eq:partiallyexpanded} until a normalized coefficient exceeds the requested error tolerance (i.e.~$b_{\max}  = \max [\mathbf{b}_{l} \odot \boldsymbol{\chi}_x^{(k)\odot-1} ]> \epsilon_x$). The last term whose coefficient does not exceed the requested error tolerance is used to define the order of the element on the ensuing mesh. As the error in Eq.~\eqref{eq:relerrorspace} is not one-to-one with the error computed from reducing the polynomial degree, in practice, it is sometimes useful to employ a small safety tolerance on the reduction step to avoid producing a mesh that dissatisfies the requested error tolerance on the next mesh iteration (i.e.~check $b_{\max} < \chi\epsilon$, where $\chi \in (0,1]$ is small). It is noted that the smallest degree to which the order of the element $k$ can be reduced is degree 1.  

\subsection{Merging Elements}

In addition to reducing the polynomial degree in each element, the mesh size can be reduced by merging adjacent elements. To determine which elements may be merged, we rely upon power series expansions of the state. Given two elements, $k$ and $k+1$, that satisfy the error tolerance $\epsilon_x$ and possess the \textit{same} number of nodal points $p^{(k)} = p^{(k+1)}$, the elements may be merged if the polynomial expansion from each element may be extended into the neighboring element and the difference between the original state approximation and the expanded polynomial is at most $\epsilon_x$. The power series expansions that will be used are expanded about the junction of the two elements, $x^{(k)}_{p^{(k)}+1} = x^{(k+1)}_1$. Expressing the Lagrange basis in terms of the power series expansions, we have 
\begin{align}
    \phi_i^{(k)}(r) &= \sum_{l = 0}^{p^{(k)}} \left(r-1\right)^la_{li}^{(k)}, \label{eq:l1} \\
    \phi_i^{(k+1)}(r)  &= \sum_{l = 0}^{p^{(k+1)}} r^la_{li}^{(k+1)}. \label{eq:l2}
\end{align}
The expansion coefficients may be determined in a similar manner to Eq.~\eqref{eq:coeffcompute}. Combining Eq.~\eqref{eq:stateapproxorig} with Eq.~\eqref{eq:l1} in element $k$ and Eq.~\eqref{eq:l2} in element $k+1$ yields
\begin{align}
    y_h^{(k)}(x) &= \sum_{l=1}^{p^{(k)}+1 }c_{l}^{(k)}\left(\frac{x-x_{p^{(k)}+1}^{(k)}}{h^{(k)}} \right)^l, \\
    c_{l}^{(k)} &= \sum_{n = 1}^{p^{(k)}+1}Y_n^{(k)}a_{ln}^{(k)}, \label{eq:exp1} \\
    y_h^{(k+1)}(x) &= \sum_{l=1}^{p^{(k+1)}+1}c_{l}^{(k+1)}\left(\frac{x-x_1^{(k+1)}}{h^{(k+1)}} \right)^l,  \\
    c_{l}^{(k+1)} &= \sum_{n = 1}^{p^{(k+1)}+1}Y_{n}^{(k+1)}a_{ln}^{(k+1)}, \label{eq:exp2}
\end{align}
where, because $p^{(k)} = p^{(k+1)}$, the equations differ only in that $h^{(k)}$ appears in the denominator of Eq.~\eqref{eq:exp1} and $h^{(k+1)}$ appears in the denominator of Eq.~\eqref{eq:exp2}. It is possible to manipulate the forms in Eqs.~\eqref{eq:exp1} and \eqref{eq:exp2} such that the denominators are identical. This can be done by 
\begin{align}
    y_h^{(k)}(x) &= \sum_{l = 1}^{p^{(k)}+1} b_{l}^{(k)}\left(\frac{x-x^{(k)}_{p^{(k)}+1}}{\bar{h}^{(k)}} \right)^l, \\
    b_{l}^{(k)} &= c_{l}^{(k)}\left(\frac{\bar{h}^{(k)}}{h^{(k)}} \right)^l, \\
    y_h^{(k+1)}(x) &= \sum_{l = 1}^{p^{(k+1)}+1} b_{l}^{(k+1)}\left(\frac{x-x_1^{(k+1)}}{\bar{h}^{(k)}} \right)^l,\\
    b_{l}^{(k+1)} &= c_{l}^{(k+1)}\left(\frac{\Delta \bar{h}^{(k)}}{h^{(k+1)}} \right)^l,
\end{align}
where $\bar{h}^{(k)} = \max{\left[h^{(k)}, h^{(k+1)}\right]}$. As the degree of both polynomials are identical, the difference is represented by 
\begin{equation}
    y_h^{(k)}(x) - y_h^{(k+1)}(x)  
    = \sum_{l = 1}^{p^{(k)}+1} \left(b_{l}^{(k)}- b_{l}^{(k+1)} \right)\left(\frac{x-x_1^{(k+1)}}{\bar{h}^{(k)}} \right)^l.
\end{equation}
Then, because $\left|x-x_1^{(k+1)}\right|/\bar{h}^{(k)} \leq 1$ for $x$ over the patch $\Omega_k\:\cup\:\Omega_{k+1}$, by the triangle inequality, we have 
\begin{equation}
    \max \left[\left|y_h^{(k)}(x) - y_h^{(k+1)}(x) \right| \right] \leq \sum_{l =1}^{p^{(k)}+1} \left|b_{l}^{(k)}- b_{l}^{(k+1)} \right|, \label{eq:absdiffstate}
\end{equation}
for $x$ over the patch. The right-hand side of Eq.~\eqref{eq:absdiffstate} may be evaluated for all points in time. To relate the error bound to the relative error tolerance, we multiply by the scaling weight $\boldsymbol{\chi}_x^{(k)\odot-1}$ for all temporal points. If the  
\begin{equation}
    \max \left[\boldsymbol{\chi}_x^{(k)\odot-1} \odot \sum_{l =1}^{p^{(k)}+1} \left|\mathbf{b}_{l}^{(k)}- \mathbf{b}_{l}^{(k+1)} \right| \right] < \epsilon_x,
\end{equation}
the two elements can be merged, where $\mathbf{b}_l^{(k)},\:\mathbf{b}_l^{(k+1)} \in \mathbb{R}^{N_t +1}$ represent the $l^{\mathrm{th}}$ coefficient in the power series expansion at every point in time over the elements $k$ and $k+1$, respectively. If at any temporal point the bound exceeds the relative error tolerance, the mesh intervals are not merged. Again, in practice, it is useful to employ a safety tolerance, $\chi \in (0,1]$, on the reduction step (i.e.~$\chi \epsilon$), as the error in Eq.~\eqref{eq:relerrorspace} is not one-to-one with the relative error produced from Eq.~\eqref{eq:absdiffstate}.

\bibliographystyle{elsarticle-num}   

\bibliography{references_AIAAJforArXiV} 



\end{document}